\documentclass[11pt]{article}
\usepackage[utf8]{inputenc}
\usepackage{geometry}

\usepackage[dvipsnames]{xcolor}

\usepackage{epsfig}
\usepackage{graphicx}
\usepackage{amsbsy}
\usepackage{amsmath}
\usepackage{amsfonts}
\usepackage{amssymb}
\usepackage{textcomp}
\usepackage{hyperref}
\usepackage{aliascnt}
\usepackage{bbm}
\usepackage{tikz-cd}
\usepackage{enumitem}\setlist[enumerate]{label=\textnormal{(\roman*)}}

\usepackage[section]{placeins}

\renewcommand{\geq}{\geqslant}

\newcommand{\mcm}[3]{\newcommand{#1}[#2]{{\ensuremath{#3}}}} 
\mcm{\tuple}{1}{\langle #1 \rangle}
\mcm{\name}{1}{\ulcorner #1 \urcorner}
\mcm{\Nbb}{0}{\mathbb{N}}
\mcm{\N}{0}{\mathbb{N}}
\mcm{\Zbb}{0}{\mathbb{Z}}
\mcm{\Rbb}{0}{\mathbb{R}}
\mcm{\Cbb}{0}{\mathbb{C}}
\mcm{\Qbb}{0}{\mathbb{Q}}
\mcm{\Sbb}{0}{\mathbb{S}}
\mcm{\Hbb}{0}{\mathbb{H}}
\mcm{\Kbb}{0}{\mathbb{K}}
\mcm{\Acal}{0}{\cal A}
\mcm{\Bcal}{0}{\mathcal B}
\mcm{\Ccal}{0}{\cal C}
\mcm{\Dcal}{0}{\cal D}
\mcm{\Ecal}{0}{\cal E}
\mcm{\Fcal}{0}{\cal F}
\mcm{\Gcal}{0}{\cal G}
\mcm{\Hcal}{0}{\cal H}
\mcm{\Ical}{0}{\cal I}
\mcm{\Jcal}{0}{\cal J}
\mcm{\Kcal}{0}{\cal K}
\mcm{\Lcal}{0}{\cal L}
\mcm{\Mcal}{0}{\cal M}
\mcm{\Ncal}{0}{\cal N}
\mcm{\Ocal}{0}{{\cal O}}
\mcm{\Pcal}{0}{{\cal P}}
\mcm{\Qcal}{0}{{\cal Q}}
\mcm{\Rcal}{0}{{\cal R}}
\mcm{\Scal}{0}{{\cal S}}
\mcm{\Tcal}{0}{{\cal T}}
\mcm{\Ucal}{0}{{\cal U}}
\mcm{\Vcal}{0}{{\cal V}}
\mcm{\Wcal}{0}{{\cal W}}
\mcm{\Xcal}{0}{{\cal X}}
\mcm{\Ycal}{0}{{\cal Y}}
\mcm{\Zcal}{0}{{\cal Z}}
\mcm{\Mfrak}{0}{\mathfrak M}

\usepackage{amsthm}

\newcommand{\theoremize}[2]{\newaliascnt{#1}{thm} \newtheorem{#1}[#1]{#2} \aliascntresetthe{#1}}

\theoremstyle{plain}
\newtheorem{thm}{Theorem}[section]
\theoremize{lem}{Lemma}

\theoremize{skolem}{Skolem}
\theoremize{fact}{Fact}
\theoremize{setting}{Setting}

\newtheorem{claim}{Claim}[thm]
\theoremize{obs}{Observation}
\theoremize{prop}{Proposition}
\theoremize{cor}{Corollary}
\theoremize{que}{Question}
\theoremize{oque}{Open Question}
\theoremize{con}{Conjecture}

\theoremstyle{definition}
\theoremize{dfn}{Definition}
\theoremize{rem}{Remark}
\theoremize{eg}{Example}
\theoremize{exercise}{Exercise}
\theoremstyle{plain}

\hypersetup{
    draft = false,
    bookmarksopen=true,
    colorlinks,
    linkcolor={red!60!black},
    citecolor={green!60!black},
    urlcolor={blue!60!black}
}

\newcommand{\sm}{\setminus}

\makeatletter
\def\thm@space@setup{\thm@preskip=2pt
\thm@postskip=2pt}
\makeatother

\makeatletter
\newenvironment{pf}[1][\proofname]{\par
  \pushQED{\qed}%
  \normalfont \topsep0\p@\relax
  \trivlist
  \item[\hskip\labelsep\itshape
  #1\@addpunct{.}]\ignorespaces
}{%
  \popQED\endtrivlist\@endpefalse
}
\makeatother

\newenvironment{cproof}{\noindent\textit{Proof of Claim.}}{\phantom{!}\!\hfill$\diamondsuit$\vspace{2pt}}

\usepackage{titlesec}

\titleformat{\section}
  {\normalfont\scshape\filcenter}{\thesection}{1em}{}

\titleformat{\subsection}
 {\normalfont\scshape\filcenter}{\thesubsection}{1em}{}

\titlespacing*{\section}
{0pt}{5pt plus 1pt minus 1pt}{5pt plus 1pt minus 1pt}

\titlespacing*{\subsection}
{0pt}{5pt plus 1pt minus 1pt}{5pt plus 1pt minus 1pt}
\usepackage{tikzit}

\tikzstyle{black dot}=[fill=black, draw=none, shape=circle, scale=0.3]
\tikzstyle{backing}=[fill=white, draw=none, shape=circle, scale=0.85]
\tikzstyle{small back}=[fill=white, draw=none, shape=circle, scale=0.95]
\tikzstyle{medium back}=[fill=white, draw=none, shape=circle, scale=1.2]

\tikzstyle{c0}=[-, line width=0.08em]
\tikzstyle{c0 arrow}=[->, line width=0.08em]
\tikzstyle{c0 bold}=[-, line width=0.15em]
\tikzstyle{c1}=[-, draw={rgb,255: red,135; green,206; blue,235}, line width=0.08em]
\tikzstyle{c1 arrow}=[->, draw={rgb,255: red,135; green,206; blue,235}, line width=0.08em]
\tikzstyle{c1 bold}=[-, draw={rgb,255: red,135; green,206; blue,235}, line width=0.15em]
\tikzstyle{c2}=[-, draw={rgb,255: red,255; green,69; blue,0}, line width=0.08em]
\tikzstyle{c2 arrow}=[->, draw={rgb,255: red,255; green,69; blue,0}, line width=0.08em]
\tikzstyle{c3}=[-, draw={rgb,255: red,102; green,176; blue,50}, line width=0.08em]
\tikzstyle{c3 arrow}=[->, draw={rgb,255: red,102; green,176; blue,50}, line width=0.08em]
\tikzstyle{d1}=[-, dashed, line width=0.08em]
\tikzstyle{d1 arrow}=[->, line width=0.08em, dashed]
\tikzstyle{c1 dashed}=[-, draw={rgb,255: red,135; green,206; blue,235}, line width=0.08em, dashed, dash pattern=on 1em off 1em]
\tikzstyle{c2 dashed}=[-, draw={rgb,255: red,255; green,69; blue,0}, line width=0.08em, dashed, dash pattern=on 1em off 1em]

\usepackage{verbatim}

\mcm{\restric}{0}{\upharpoonright}

\mcm{\minorp}{0}{\ \preccurlyeq \ }
\DeclareRobustCommand{\minorm}{\text{\reflectbox{$\ \preccurlyeq \ $}}}

\newcommand{\rcol}{\textcolor{BrickRed}{\textsc{red}}}
\newcommand{\bcol}{\textcolor{SkyBlue}{\textsc{blue}}}
\newcommand{\gcol}{\textcolor{Green}{\textsc{green}}}

\newcommand{\se}{\subseteq}

\mcm{\Fbb}{0}{\mathbb{F}}

\title{\vspace{-1cm}Hardness of Planarity for Weak Temporal Sequences of 2-Connected Graphs}
\author{Johannes Carmesin\thanks{TU Freiberg, funded by DFG, project number 546892829.} \ and Will J.\ Turner\footnotemark[1]}
\date{}

\begin{document}

\maketitle

\vspace{-0.75cm}
\begin{center}
\begin{minipage}{0.9\textwidth}
\small
\textbf{Abstract.} A \emph{weak deletion sequence} is a sequence $(G_1,\ldots,G_n)$ of graphs so that for each $i\in[n-1]$ either $G_i$ is isomorphic to a subgraph of $G_{i+1}$, or vice versa: $G_{i+1}$ is isomorphic to a subgraph of $G_i$. We prove that determining the simultaneous planar embeddability of weak deletion sequences of $2$-connected graphs is NP-hard. 


\end{minipage}
\end{center}

\section{Introduction}

This paper forms part of the follow-up work to \cite{temporalsequences2025} along with \cite{fpttreesofgraphs2025} and \cite{contractionsequences2025} which study temporal sequences, a recent generalisation of temporal graphs; we explain this connection in the concluding remarks. 
A \emph{temporal (deletion) sequence} is a sequence $(G_1,\ldots,G_n)$ of graphs so that for each $i\in[n-1]$ the graph $G_i$ is obtained from $G_{i+1}$ by either deleting edges and isolated vertices or adding edges and isolated vertices, but not both. There are two natural variants: in \emph{strict temporal sequences} no relabelling is permitted, while in \emph{weak temporal sequences} graphs are considered up to isomorphism. This paper is concerned with weak temporal sequences. 

Informally, a temporal sequence $(G_1,\dots,G_n)$ is simultaneously embeddable if one can choose plane embeddings of the graphs $G_i$ that are compatible along the sequence.
Formally, a weak temporal sequence $(G_1,\dots,G_n)$ is said to be \emph{weakly simultaneously embeddable} if there exists a sequence of planar embeddings $(\iota_1,\dots,\iota_n)$ such that for each $i \in [n-1]$, whenever one of $G_i,G_{i+1}$ is isomorphic to a subgraph of the other, the embedding of the larger graph induces, up to isomorphism, the embedding of the smaller.
The main result of this paper is the following. 

\begin{prop}\label{thm:main}
    The problem of deciding whether weak deletion sequences of $2$-connected graphs admit weak simultaneous embeddings is NP-hard.
\end{prop}

Perhaps surprisingly, the distinction between weak and strict temporal sequences has striking algorithmic consequences. In contrast to \autoref{thm:main}, the natural strict analogue---deciding simultaneous embeddability for strict deletion sequences of $2$-connected graphs---is solvable in polynomial time~\cite{temporalsequences2025}.

In the landscape of tractability, it is the $2$-connected regime that exhibits the richest structure. For sequences of $3$-connected graphs, the problem is trivial (and hence in~P), whereas for connected graphs in general it is NP-hard---both for the weak and the strict variants---by the result of~\cite{gassner2006simultaneous}.

The remainder of the paper is structured as follows. 
The proof of \autoref{thm:main} will be found in \autoref{sec:part1} and the supporting lemmas are subsequently proved in \autoref{sec:part2} and \autoref{sec:part3}. In \autoref{sec:conc}, we make some concluding remarks. In the appendix, we will adapt the proof of \autoref{thm:main} to prove NP-hardness for a different weakening of strict temporal sequences called `indefinite temporal sequences', where we forgo the assumption that for each edge to which a minor operation is applied, we fix in advance what that minor operation is going to be.

\section{Hybrid temporal sequences}
In this section we reduce our problem to one concerning more general objects.
Reasoning about graph isomorphisms can get very messy and, for this reason, we will instead work with a `hybrid' form of deletion sequence defined below. First however, we exhibit the following simple example which highlights the difference between strict and weak temporal sequences. 

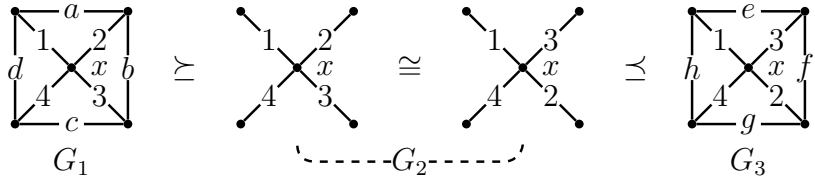
\begin{figure}
    \centering
        { \large\scalebox{1}{\begin{tikzpicture}
	\begin{pgfonlayer}{nodelayer}
		\node [style=none] (420) at (-2.5, -6) {$G_1$};
		\node [style=black dot] (427) at (-4, -2) {};
		\node [style=black dot] (428) at (-1, -2) {};
		\node [style=black dot] (429) at (-4, -5) {};
		\node [style=black dot] (430) at (-1, -5) {};
		\node [style=black dot] (431) at (-2.5, -3.5) {};
		\node [style=none] (432) at (0.5, -3.5) {$\succeq$};
		\node [style=none] (439) at (-1.75, -3.5) {$x$};
		\node [style=backing] (452) at (-2.5, -2) {};
		\node [style=none] (455) at (-2.5, -2) {$a$};
		\node [style=backing] (456) at (-1, -3.5) {};
		\node [style=none] (457) at (-1, -3.5) {$b$};
		\node [style=backing] (458) at (-2.5, -5) {};
		\node [style=none] (459) at (-2.5, -5) {$c$};
		\node [style=backing] (460) at (-4, -3.5) {};
		\node [style=none] (461) at (-4, -3.5) {$d$};
		\node [style=backing] (462) at (-3.25, -2.75) {};
		\node [style=none] (463) at (-3.25, -2.75) {$1$};
		\node [style=backing] (464) at (-1.75, -2.75) {};
		\node [style=none] (465) at (-1.75, -2.75) {$2$};
		\node [style=backing] (468) at (-3.25, -4.25) {};
		\node [style=none] (469) at (-3.25, -4.25) {$4$};
		\node [style=backing] (470) at (-1.75, -4.25) {};
		\node [style=none] (471) at (-1.75, -4.25) {$3$};
		\node [style=none] (472) at (6.5, -6) {$G_2$};
		\node [style=black dot] (473) at (2, -2) {};
		\node [style=black dot] (474) at (5, -2) {};
		\node [style=black dot] (475) at (2, -5) {};
		\node [style=black dot] (476) at (5, -5) {};
		\node [style=black dot] (477) at (3.5, -3.5) {};
		\node [style=none] (478) at (4.25, -3.5) {$x$};
		\node [style=backing] (487) at (2.75, -2.75) {};
		\node [style=none] (488) at (2.75, -2.75) {$1$};
		\node [style=backing] (489) at (4.25, -2.75) {};
		\node [style=none] (490) at (4.25, -2.75) {$2$};
		\node [style=backing] (491) at (2.75, -4.25) {};
		\node [style=none] (492) at (2.75, -4.25) {$4$};
		\node [style=backing] (493) at (4.25, -4.25) {};
		\node [style=none] (494) at (4.25, -4.25) {$3$};
		\node [style=none] (495) at (6.5, -3.5) {$\cong$};
		\node [style=black dot] (497) at (8, -2) {};
		\node [style=black dot] (498) at (11, -2) {};
		\node [style=black dot] (499) at (8, -5) {};
		\node [style=black dot] (500) at (11, -5) {};
		\node [style=black dot] (501) at (9.5, -3.5) {};
		\node [style=none] (502) at (10.25, -3.5) {$x$};
		\node [style=backing] (503) at (8.75, -2.75) {};
		\node [style=none] (504) at (8.75, -2.75) {$1$};
		\node [style=backing] (505) at (10.25, -2.75) {};
		\node [style=none] (506) at (10.25, -2.75) {$3$};
		\node [style=backing] (507) at (8.75, -4.25) {};
		\node [style=none] (508) at (8.75, -4.25) {$4$};
		\node [style=backing] (509) at (10.25, -4.25) {};
		\node [style=none] (510) at (10.25, -4.25) {$2$};
		\node [style=none] (511) at (15.5, -6) {$G_3$};
		\node [style=black dot] (512) at (14, -2) {};
		\node [style=black dot] (513) at (17, -2) {};
		\node [style=black dot] (514) at (14, -5) {};
		\node [style=black dot] (515) at (17, -5) {};
		\node [style=black dot] (516) at (15.5, -3.5) {};
		\node [style=none] (517) at (12.5, -3.5) {$\preceq$};
		\node [style=none] (518) at (16.25, -3.5) {$x$};
		\node [style=backing] (519) at (15.5, -2) {};
		\node [style=none] (520) at (15.5, -2) {$e$};
		\node [style=backing] (521) at (17, -3.5) {};
		\node [style=none] (522) at (17, -3.5) {$f$};
		\node [style=backing] (523) at (15.5, -5) {};
		\node [style=none] (524) at (15.5, -5) {$g$};
		\node [style=backing] (525) at (14, -3.5) {};
		\node [style=none] (526) at (14, -3.5) {$h$};
		\node [style=backing] (527) at (14.75, -2.75) {};
		\node [style=none] (528) at (14.75, -2.75) {$1$};
		\node [style=backing] (529) at (16.25, -2.75) {};
		\node [style=none] (530) at (16.25, -2.75) {$3$};
		\node [style=backing] (531) at (14.75, -4.25) {};
		\node [style=none] (532) at (14.75, -4.25) {$4$};
		\node [style=backing] (533) at (16.25, -4.25) {};
		\node [style=none] (534) at (16.25, -4.25) {$2$};
		\node [style=none] (535) at (6, -6) {};
		\node [style=none] (536) at (7, -6) {};
		\node [style=none] (537) at (4, -6) {};
		\node [style=none] (538) at (3.5, -5.5) {};
		\node [style=none] (539) at (9, -6) {};
		\node [style=none] (540) at (9.5, -5.5) {};
	\end{pgfonlayer}
	\begin{pgfonlayer}{edgelayer}
		\draw [style=c0] (427) to (431);
		\draw [style=c0] (431) to (428);
		\draw [style=c0] (428) to (427);
		\draw [style=c0] (427) to (429);
		\draw [style=c0] (429) to (431);
		\draw [style=c0] (430) to (428);
		\draw [style=c0] (431) to (430);
		\draw [style=c0] (430) to (429);
		\draw [style=c0] (473) to (477);
		\draw [style=c0] (477) to (474);
		\draw [style=c0] (475) to (477);
		\draw [style=c0] (477) to (476);
		\draw [style=c0] (497) to (501);
		\draw [style=c0] (501) to (498);
		\draw [style=c0] (499) to (501);
		\draw [style=c0] (501) to (500);
		\draw [style=c0] (512) to (516);
		\draw [style=c0] (516) to (513);
		\draw [style=c0] (513) to (512);
		\draw [style=c0] (512) to (514);
		\draw [style=c0] (514) to (516);
		\draw [style=c0] (515) to (513);
		\draw [style=c0] (516) to (515);
		\draw [style=c0] (515) to (514);
		\draw [style=d1] (535.center)
			 to (537.center)
			 to [bend right=315, looseness=1.25] (538.center);
		\draw [style=d1] (536.center)
			 to (539.center)
			 to [bend right=45] (540.center);
	\end{pgfonlayer}
\end{tikzpicture}}}

        \caption{A strict temporal sequence $(G_1,G_2,G_3)$. The graph $G_2$ is obtained from $G_1$ by deleting the edges $a,b,c$ and $d$, and obtained from $G_3$ by deleting the edges $e,f,g$ and $h$. Also depicted is an embedding of $G_1$, an embedding of $G_3$, and the respective embeddings that these induce on $G_2$.}
    \label{fig:iso_easo}
\end{figure}

\begin{eg}
    Consider the strict temporal sequence $(G_1,G_2,G_3)$ depicted in \autoref{fig:iso_easo}. As a strict temporal sequence, it admits no simultaneous embedding. Indeed, each of the possible embeddings for $G_1$ induces the cycling order $(1234)$ or its inverse on the edges incident with the vertex $x$. However, each of the possible embeddings for $G_3$ induces the cyclic order $(1243)$ or its inverse on the edges incident with the vertex $x$. Hence no embeddings of $G_1$ and $G_3$ can ever induce the same embedding of $G_2$ and we conclude that $(G_1,G_2,G_3)$ admits no (strict) simultaneous embedding. In contrast, if we instead view $(G_1,G_2,G_3)$ as a weak temporal sequence, then we find that it admits a weak simultaneous embedding. Indeed, we can arbitrarily permute the cyclic order of the edges incident with $x$ in embeddings of $G_2$ by considering appropriate automorphisms.
\end{eg}

\begin{dfn}[Hybrid deletion sequence]
    Let $H$ and $G$ be graph and $A\subseteq E(H)$. We say that $H$ is a \emph{hybrid subgraph of $G$ w.r.t.\ $A$} if $H$ is obtained from a subgraph of $G$ by relabeling vertices and relabeling some edges outside of $A$. We call $A$ the set of \emph{strict edges} of $H$ (resp.\ $G$). Observe that $H$ is a (strict) subgraph of $G$ exactly when $A=E(H)$.

    A \emph{hybrid deletion sequence} is a temporal sequence for the hybrid subgraph relation. Let $(G_1,\ldots,G_n)$ be a hybrid deletion sequence. The union\footnote{Here we adhere to the convention that each edge label that appears in our sequence is either strict in all hybrid subgraph relations or not strict in all the hybrid subgraph relations.} of the sets of strict edges for each~$G_i$ is called the set of \emph{strict edges} of $(G_1,\ldots,G_n)$.

    For a graph isomorphism $\theta:G\rightarrow H$, we say that $\theta$ \emph{fixes} some edge set $A\subseteq E(G)$ if $\theta$ is the identity on the edge labels $A$ and preserves the incident relations between the edges $A$. Now let $(G_1,\ldots,G_n)$ be a hybrid deletion sequence with strict edges $A$ and let $(\phi_1,\ldots,\phi_n)$ be a respective sequence of embeddings. We say that $(\phi_1,\ldots,\phi_n)$ is a \emph{(hybrid) simultaneous embedding} of $(G_1,\ldots,G_n)$ if for each $i\in[n-1]$ whenever $G_i$ is a hybrid subgraph of $G_{i+1}$ for some $A\subseteq E(G_i)$, we have that $\phi_i$ is the same up to $A$-fixing isomorphism as the embedding induced by $\phi_{i+1}$ on some subgraph of $G_{i+1}$; and otherwise $G_{i+1}$ is a hybrid subgraph of $G_i$ for some $A\subseteq E(G_{i+1})$, in which case we have that $\phi_{i+1}$ is the same up to $A$-fixing isomorphism as the embedding induced by $\phi_{i}$ on some subgraph of $G_{i}$.
\end{dfn}
A hybrid deletion sequence has fewer possible simultaneous embeddings in general than a weak deletion sequence containing the same graphs.

\begin{lem}\label{lem:hybrid_equivalence}
    For every hybrid deletion sequence $\Gcal=(G_1,\ldots,G_n)$ of $2$-connected graphs, there exists a weak deletion sequence $\Gcal'=(G_1',\ldots,G_n')$ of $2$-connected graphs so that $\Gcal$ admits a hybrid simultaneous embedding if and only if $\Gcal'$ admits a weak simultaneous embedding. Furthermore, $\text{size}(\Gcal')$ is bounded by a polynomial in terms of $\text{size}(\Gcal)$.
\end{lem}
\begin{pf}[Sketch of proof]
    Let $A$ denote the set of strict edges in $\Gcal$. For each $i\in[n]$, we obtain $G_i'$ from $G_i$ by replacing each edge $e$ in $G_i$ which is not strict by a `gadget' $Q_e$. The gadgets $Q_e$ are planar graphs which are chosen so that they are pairwise not subgraphs of one another, so that the graphs $G_i$ remain $2$-connected and so that the size of each $Q_e$ is bounded by a polynomial in terms of the number of edges in $\Gcal$. 
    In this way, we ensure that whenever $G_i$ is a hybrid subgraph of $G_{i+1}$ for $A\subseteq E(G_i)$, we get that the only isomorphisms from $G_{i}'$ to a minor of $G_{i+1}'$ are those which fix the image of $A$ under our replacements.

    It remains to give a construction for the graphs $Q_e$ that satisfies the above criteria. For example, we could take a graph $Q_e$ to be the graph obtained from a grid, whose size is a constant factor larger than the size of $\Gcal$, by adding a few additional edges so as to make each $Q_e$ distinct. When we replace $e$ by $Q_e$, we attach the two opposite corners of the grid to the endvertices of $e$.
\end{pf}

The proof of \autoref{thm:main} is divided into three sections. In \autoref{sec:part1}, we show that if three classes of hybrid deletion sequences with certain properties exist, then $3$-SAT can be reduced to the simultaneous embedding problem for $2$-connected hybrid deletion sequences. In Sections \ref{sec:part2} and~\ref{sec:part3}, we construct the required three classes: In \autoref{sec:part2} we first construct these classes in a very special case, and then in \autoref{sec:part3} we show how these special cases can be combined to give the complete classes. Finally, in \autoref{sec:conc} we remark upon further directions for study.
For general graph-theoretic and graph-embedding terminology, we follow the book of Mohar and Thomassen~\cite{graphsonsurfaces2001}. For terminology specific to temporal sequences and simultaneous embeddings, we refer the reader to~\cite{temporalsequences2025}.

\section{Reducing $3$-SAT to a construction problem for hybrid deletion sequences}\label{sec:part1}

In this section, we reduce $3$-SAT to deciding hybrid simultaneous embeddability, given that it is possible to construct certain hybrid deletion sequences.

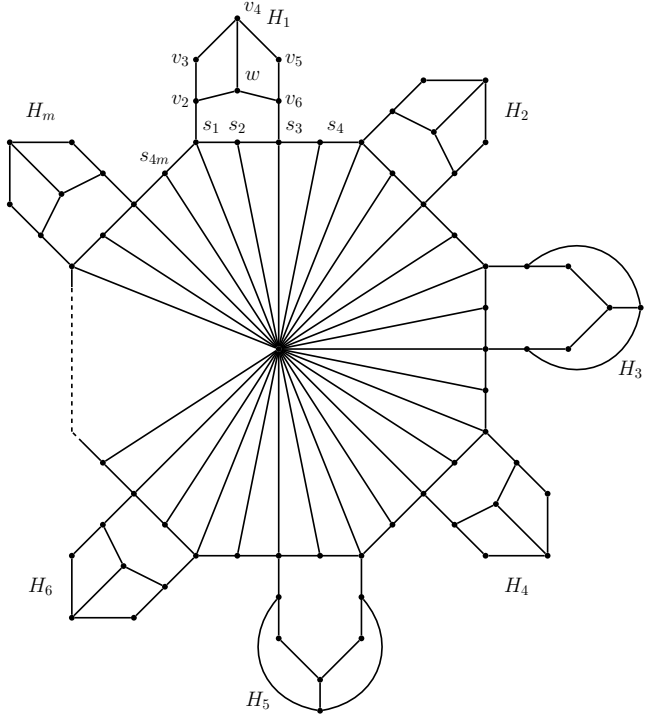
\begin{figure}
    \begin{minipage}{0.59\textwidth}
        \Large\scalebox{0.55}{\begin{tikzpicture}
	\begin{pgfonlayer}{nodelayer}
		\node [style=black dot] (0) at (-3, 2) {};
		\node [style=black dot] (1) at (-1, 2) {};
		\node [style=black dot] (2) at (1, 2) {};
		\node [style=black dot] (3) at (3, 2) {};
		\node [style=black dot] (4) at (5, 2) {};
		\node [style=black dot] (5) at (6.5, 0.5) {};
		\node [style=black dot] (6) at (8, -1) {};
		\node [style=black dot] (7) at (9.5, -2.5) {};
		\node [style=black dot] (8) at (11, -4) {};
		\node [style=black dot] (9) at (11, -6) {};
		\node [style=black dot] (10) at (11, -8) {};
		\node [style=black dot] (11) at (11, -10) {};
		\node [style=black dot] (12) at (8, -15) {};
		\node [style=black dot] (13) at (6.5, -16.5) {};
		\node [style=black dot] (14) at (5, -18) {};
		\node [style=black dot] (15) at (3, -18) {};
		\node [style=black dot] (16) at (1, -18) {};
		\node [style=black dot] (17) at (-1, -18) {};
		\node [style=black dot] (18) at (-3, -18) {};
		\node [style=black dot] (22) at (-4.5, 0.5) {};
		\node [style=black dot] (23) at (-6, -1) {};
		\node [style=black dot] (24) at (-9, -4) {};
		\node [style=black dot] (25) at (11, -12) {};
		\node [style=black dot] (26) at (-7.5, -2.5) {};
		\node [style=black dot] (27) at (9.5, -13.5) {};
		\node [style=black dot] (32) at (-4.5, -16.5) {};
		\node [style=black dot] (33) at (-6, -15) {};
		\node [style=black dot] (34) at (-7.5, -13.5) {};
		\node [style=black dot] (35) at (1, -8) {};
		\node [style=none] (39) at (-9, -12) {};
		\node [style=none] (40) at (-9, -5) {};
		\node [style=black dot] (41) at (-3, 4) {};
		\node [style=black dot] (42) at (-3, 6) {};
		\node [style=black dot] (43) at (1, 4) {};
		\node [style=black dot] (44) at (1, 6) {};
		\node [style=black dot] (45) at (-1, 8) {};
		\node [style=black dot] (46) at (-1, 4.5) {};
		\node [style=black dot] (47) at (6.5, 3.5) {};
		\node [style=black dot] (48) at (8, 5) {};
		\node [style=black dot] (49) at (9.5, 0.5) {};
		\node [style=black dot] (50) at (11, 2) {};
		\node [style=black dot] (51) at (11, 5) {};
		\node [style=black dot] (52) at (8.5, 2.5) {};
		\node [style=black dot] (53) at (15, -4) {};
		\node [style=black dot] (54) at (13, -4) {};
		\node [style=black dot] (55) at (13, -8) {};
		\node [style=black dot] (56) at (15, -8) {};
		\node [style=black dot] (57) at (17, -6) {};
		\node [style=black dot] (58) at (18.5, -6) {};
		\node [style=black dot] (59) at (12.5, -13.5) {};
		\node [style=black dot] (60) at (14, -15) {};
		\node [style=black dot] (61) at (9.5, -16.5) {};
		\node [style=black dot] (62) at (11, -18) {};
		\node [style=black dot] (63) at (14, -18) {};
		\node [style=black dot] (64) at (11.5, -15.5) {};
		\node [style=black dot] (65) at (-10.5, -2.5) {};
		\node [style=black dot] (66) at (-7.5, 0.5) {};
		\node [style=black dot] (67) at (-12, -1) {};
		\node [style=black dot] (68) at (-9, 2) {};
		\node [style=black dot] (69) at (-12, 2) {};
		\node [style=black dot] (70) at (-9.5, -0.5) {};
		\node [style=black dot] (71) at (-7.5, -16.5) {};
		\node [style=black dot] (72) at (-9, -18) {};
		\node [style=black dot] (73) at (-4.5, -19.5) {};
		\node [style=black dot] (74) at (-6, -21) {};
		\node [style=black dot] (75) at (-9, -21) {};
		\node [style=black dot] (76) at (-6.5, -18.5) {};
		\node [style=black dot] (77) at (1, -20) {};
		\node [style=black dot] (78) at (5, -20) {};
		\node [style=black dot] (79) at (5, -22) {};
		\node [style=black dot] (80) at (1, -22) {};
		\node [style=black dot] (81) at (3, -24) {};
		\node [style=black dot] (82) at (3, -25.5) {};
		\node [style=none] (83) at (-8.5, -12.5) {};
		\node [style=none] (84) at (-2.25, 2.75) {$s_1$};
		\node [style=none] (85) at (-1, 2.75) {$s_2$};
		\node [style=none] (86) at (1.75, 2.75) {$s_3$};
		\node [style=none] (87) at (3.75, 2.75) {$s_4$};
		\node [style=none] (88) at (-5, 1.25) {$s_{4m}$};
		\node [style=none] (89) at (-3.75, 4) {$v_2$};
		\node [style=none] (90) at (-3.75, 6) {$v_3$};
		\node [style=none] (91) at (-0.25, 8.5) {$v_4$};
		\node [style=none] (92) at (1.75, 6) {$v_5$};
		\node [style=none] (93) at (1.75, 4) {$v_6$};
		\node [style=none] (94) at (-0.25, 5.25) {$w$};
		\node [style=none] (95) at (1, 8) {$H_1$};
		\node [style=none] (96) at (12.5, 3.5) {$H_2$};
		\node [style=none] (97) at (18, -9) {$H_3$};
		\node [style=none] (98) at (12.5, -19.5) {$H_4$};
		\node [style=none] (99) at (0, -25) {$H_5$};
		\node [style=none] (100) at (-10.5, -19.5) {$H_6$};
		\node [style=none] (101) at (-10.5, 3.5) {$H_m$};
	\end{pgfonlayer}
	\begin{pgfonlayer}{edgelayer}
		\draw [style=c0] (35) to (0);
		\draw [style=c0] (1) to (35);
		\draw [style=c0] (35) to (2);
		\draw [style=c0] (35) to (3);
		\draw [style=c0] (35) to (4);
		\draw [style=c0] (35) to (5);
		\draw [style=c0] (6) to (35);
		\draw [style=c0] (35) to (7);
		\draw [style=c0] (8) to (35);
		\draw [style=c0] (35) to (9);
		\draw [style=c0] (10) to (35);
		\draw [style=c0] (35) to (11);
		\draw [style=c0] (25) to (35);
		\draw [style=c0] (35) to (27);
		\draw [style=c0] (12) to (35);
		\draw [style=c0] (0) to (4);
		\draw [style=c0] (4) to (8);
		\draw [style=c0] (8) to (25);
		\draw [style=c0] (25) to (14);
		\draw [style=c0] (14) to (18);
		\draw [style=c0] (24) to (0);
		\draw [style=c0] (22) to (35);
		\draw [style=c0] (35) to (23);
		\draw [style=c0] (26) to (35);
		\draw [style=c0] (24) to (35);
		\draw [style=c0] (24) to (40.center);
		\draw [style=d1] (39.center) to (40.center);
		\draw [style=c0] (35) to (34);
		\draw [style=c0] (33) to (35);
		\draw [style=c0] (35) to (32);
		\draw [style=c0] (18) to (35);
		\draw [style=c0] (35) to (17);
		\draw [style=c0] (16) to (35);
		\draw [style=c0] (14) to (35);
		\draw [style=c0] (35) to (15);
		\draw [style=c0] (13) to (35);
		\draw [style=c0] (24) to (67);
		\draw [style=c0] (67) to (69);
		\draw [style=c0] (69) to (68);
		\draw [style=c0] (68) to (23);
		\draw [style=c0] (0) to (42);
		\draw [style=c0] (42) to (45);
		\draw [style=c0] (45) to (44);
		\draw [style=c0] (44) to (2);
		\draw [style=c0] (4) to (48);
		\draw [style=c0] (48) to (51);
		\draw [style=c0] (51) to (50);
		\draw [style=c0] (50) to (6);
		\draw [style=c0] (8) to (53);
		\draw [style=c0] (53) to (57);
		\draw [style=c0] (57) to (56);
		\draw [style=c0] (56) to (10);
		\draw [style=c0] (25) to (60);
		\draw [style=c0] (60) to (63);
		\draw [style=c0] (63) to (62);
		\draw [style=c0] (62) to (12);
		\draw [style=c0] (14) to (79);
		\draw [style=c0] (79) to (81);
		\draw [style=c0] (81) to (80);
		\draw [style=c0] (80) to (16);
		\draw [style=c0] (33) to (71);
		\draw [style=c0] (71) to (72);
		\draw [style=c0] (72) to (75);
		\draw [style=c0] (74) to (73);
		\draw [style=c0] (74) to (75);
		\draw [style=c0] (73) to (18);
		\draw [style=c0] (69) to (70);
		\draw [style=c0] (70) to (65);
		\draw [style=c0] (70) to (66);
		\draw [style=c0] (45) to (46);
		\draw [style=c0] (46) to (41);
		\draw [style=c0] (46) to (43);
		\draw [style=c0] (47) to (52);
		\draw [style=c0] (52) to (49);
		\draw [style=c0] (52) to (51);
		\draw [style=c0, bend left=60, looseness=1.25] (54) to (58);
		\draw [style=c0, bend left=60, looseness=1.25] (58) to (55);
		\draw [style=c0] (58) to (57);
		\draw [style=c0] (64) to (59);
		\draw [style=c0] (64) to (61);
		\draw [style=c0] (64) to (63);
		\draw [style=c0, bend left=60, looseness=1.25] (82) to (77);
		\draw [style=c0, bend right=60, looseness=1.25] (82) to (78);
		\draw [style=c0] (82) to (81);
		\draw [style=c0] (71) to (76);
		\draw [style=c0] (76) to (73);
		\draw [style=c0] (76) to (75);
		\draw [style=c0] (18) to (34);
		\draw [style=c0] (34) to (83.center);
		\draw [style=d1] (83.center) to (39.center);
	\end{pgfonlayer}
\end{tikzpicture}}
    \end{minipage}
    \hfill%
    \begin{minipage}{0.39\textwidth}
        \caption{A village with $m$ houses. Observe that for the illustrated embedding $\sigma$, the inhabitants of $H_1$ and $H_2$ are `embedded inside' of their houses and the inhabitant of $H_3$ is `embedded outside'. Hence $T_\sigma(H_1)=1,T_\sigma(H_2)=1$ and $T_\sigma(H_3)=0$.}
    \label{fig:village}
    \end{minipage}
\end{figure}

A \lq village\rq\ is a particular planar graph, which is depicted in \autoref{fig:village}, and which we will formally define in the following.

A \emph{house} is the graph obtained from a path of edge-length $6$ with vertices $v_1,\ldots,v_7$ in that order by gluing the leaves of a claw to $v_2,v_4,v_6$. The central vertex of the claw is labeled $w$.
We refer to the glued claw as the \emph{inhabitant} of the house.
Given $m\in \Nbb$ and $I\se [m]$, a \emph{village} with \emph{index set} $I$ is obtained from a wheel of rim-length $4m$ as follows. Denote the rim vertices by $s_1,\ldots,s_{4m}$.
For each $i\in I$, we glue a house with its vertex $v_1$ identified to the rim vertex $s_{4i-3}$ and its vertex $v_7$ identified with the surface vertex $s_{4i-1}$. By $V(m,I)$ we denote the village with index set $I$ on a planet of equator $m$, and its $i$-th house by $H_i$. We refer to the subgraph corresponding to the original wheel as the \emph{planet (of equator $4m$)} and its rim as the \emph{surface of the planet}.

For each house $H_i$, the path $v_1\ldots v_7$ forms a cycle with the path $s_{4i-3}\ldots s_{4i-1}$ in the surface of the planet. We call this cycle the \emph{bounding cycle} for the house $H_i$. In an embedding $\sigma$ of $V(m,I)$, the bounding cycle $C_i$ for $H_i$ partitions $\Sbb^2$ into two topological disks by the Jordan Curve Theorem. In the context of $\sigma$, we say that the disk containing the centre of the planet is the \emph{outside} of $C_i$ and the other disk is the \emph{inside}.
In $\sigma$, the inhabitant of $H_i$ is embedded either on the inside or outside of the bounding cycle and, in these cases, we say the inhabitant is \emph{outside} or \emph{inside} of its house, respectively.
Given $i\in I$, we set $T_\sigma(H_i)=1$ if the inhabitant of $H_i$ is embedded inside $H_i$ in $\sigma$, and otherwise we set $T_\sigma(H_i)=0$; see \autoref{fig:village}. 

\begin{eg}
Given $m\in \Nbb$ and $I\se [m]$, for any function $f$ from $I$ to $\{0,1\}$, there is an embedding $\sigma$ of $V(m,I)$ so that for all $i\in I$ we have $T_\sigma(H_i)=f(i)$. 
\end{eg}

Given $m\in \Nbb$ and $I\se [m]$, a \emph{housing sequence} with housing index set $I$ on a planet of equator $m$ is a 2-connected hybrid deletion sequence $\Gcal$ that begins and ends with graphs isomorphic to $V(m,I)$ having the property that, in all its simultaneous embeddings, the first and last graph have the same embedding. Furthermore, the planet of $V(m,I)$ is a subgraph of all the graphs in the sequence, and every vertex of a graph in the sequence that is not equal to the centre has strictly smaller degree than the centre. 
The \emph{houses} of a housing sequence are the houses of its first graph and, given a simultaneous embedding $\Sigma=(\sigma_1,\ldots,\sigma_n)$ and an index $i\in I$, we write $T_\Sigma(H_i)$ for the truth value of the $i$-th house $H_i$ of the first graph; in formulas: $T_\Sigma(H_i):=T_{\sigma_1}(H_i)$.

\begin{lem}\label{sat-prep}
For every $3$-SAT formula $F$ of size $f$ there is a housing sequence $\Gcal=(V(m,[m]),\break G_1,\ldots,G_n, V(m,[m]))$ with equator $m=2 f$ of length at most $95$ and size at most $55000m^{18}$ so that $F$ is satisfiable if and only if $\Gcal$ admits a simultaneous embedding.
Moreover, $\Gcal$ can be constructed from $F$ in polynomial time. 
\end{lem}

We refer to the problem of determining for 2-connected hybrid deletion sequences of length at most $95$ whether they admit  simultaneous embeddings as $HDS(2,95)$.

\begin{prop}\label{HARDNESS1}
The problem $HDS(2,95)$ is NP-hard.
\end{prop}

\begin{pf}[Proof that \autoref{sat-prep} implies \autoref{HARDNESS1}.]
It suffices to prove that $3$-SAT can be reduced to $HDS(2,95)$. So let a $3$-SAT formula $F$ be given and denote its length by $f$.
By \autoref{sat-prep} we can construct in polynomial time a housing sequence $\Gcal$ with $m=2\cdot f$ houses of length at most $95$ and size at most $55000m^{18}$
so that $F$ is satisfiable if and only if $\Gcal$ admits a simultaneous embedding.
This completes the description of a reduction of $3$-SAT to $HDS(2,95)$.
\end{pf}

Clearly, combining \autoref{HARDNESS1} with \autoref{lem:hybrid_equivalence} implies \autoref{thm:main}.

The rest of this paper is dedicated to the proof of \autoref{sat-prep}. For this, we first state a few lemmas and show how they imply \autoref{sat-prep}, and then we continue by proving these lemmas. 
Given a housing sequence $\Gcal$ with housing index set $I$, the \emph{allocation function} of a simultaneous embedding $\Sigma$ of $\Gcal$ is the function $f$ from $I$ to $\{0,1\}$ defined via $f(i)=T_\Sigma(H_i)$. 
The \emph{allocation set} of $\Gcal$ is the set of all allocation functions for simultaneous embeddings of $\Gcal$. 
We shall abbreviate \lq housing sequence with housing index $I$ and equator $m$\rq\ by \lq $(I,m)$-housing sequence\rq\ (read: I'm housing sequence).

\begin{lem}\label{joint-allocation}
Given two $(I,m)$-housing sequences $\Gcal_1$ and $\Gcal_2$ with allocation sets $A_1$ and $A_2$, their concatenation $\Gcal_1\Gcal_2$ is an $(I,m)$-housing sequence with allocation set $A_1\cap A_2$. \qed
\end{lem}

In the context of an index set $I$ and a set $P\se \binom{I}{2}$ of (unordered) pairs,
define $EQ(P)$ (read: equip) to consist of those functions $f$ from $I$ to $\{0,1\}$ satisfying $f(i)=f(j)$ for all $(i,j)\in P$. 

\begin{eg}
Define an equivalence relation on $I$ as the transitive closure of the relation given by $i\sim j$ if $(i,j)\in P$. Then $EQ(P)$ consists of those functions  $f$ from $I$ to $\{0,1\}$ that are constant on each equivalence class. 
\end{eg}

\begin{lem}\label{EQ-statement}
For every $m\in \Nbb$ and $P\se \binom{[m]}{2}$ and $s\geq m^3$, there is an $([m],s)$-housing sequence of length at most $27$ and size at most $400\,s^{4/3}$ whose allocation set is $EQ(P)$. 
\end{lem}

In the context of an index set $I$ and a set $J\se \binom{I}{2}$ of (unordered) pairs,
define $NEQ(J)$ (read: not equal for $J$) to consist of those function $f$ from $I$ to $\{0,1\}$ satisfying $f(i)\neq f(j)$ for all $(i,j)\in J$. 

\begin{lem}\label{NEQ-statement}
For every $m\in \Nbb$ and $J\se \binom{[m]}{2}$ and $s\geq m^3$ there is an $([m],s)$-housing sequence of length at most $27$ and size at most $400\,s^{4/3}$ whose allocation set is $NEQ(J)$. 
\end{lem}

In the context of an index set $I$ and a set $K\se I^3$ of (ordered) triples of elements from $I$,
define $OR(K)$ (read: orc) to consist of those functions $f$ from $I$ to $\{0,1\}$ satisfying $f(i)\vee f(j)\vee f(k)=1$ for all $(i,j,k)\in K$. 
In the following, we set $K:=\{(3i-2,3i-1,3i)| \ i \text{\ divides\ } m\}$. 

\begin{lem}\label{OR-statement}
For every $m\in \Nbb$, $s\geq 27m^9$ and $K\se [m]^3$, there is an $([m],s)$-housing sequence of length at most $43$ and size at most $2001s^2$ whose allocation set is $OR(K)$.
\end{lem}

%
%

\begin{pf}[Proof of \autoref{sat-prep} assuming \autoref{EQ-statement}, \autoref{NEQ-statement} and \autoref{OR-statement}]
Let a $3$-SAT formula $F$ of size $m$ be given. Recall that its size is the number of its literals counted with multiplicity. Let $s=27m^3$.
Let $(X_i|i\in [m])$ be the family of literals of $F$. 
Next we define $([m],s)$-housing sequences $\Gcal_1,\Gcal_2$ and $\Gcal_3$ using the lemmas stated above.

{\bf Definition of $\Gcal_1$.}
We let $P$ be the set of pairs $(i,j)$ with $i\neq j$ so that $X_i=X_j$. 
By \autoref{EQ-statement} there is an $([m],s)$-housing sequence of length at most $27$ and size at most $400s^{4/3}$ whose allocation set is $EQ_s(P)$. Denote such a housing sequence by $\Gcal_1$.

{\bf Definition of $\Gcal_2$.}
For every boolean variable $X$ of $F$ so that both $X$ and its negation appear in $F$, we choose a pair $(i,j)$ so that one of $X_i$ and $X_j$ is equal to $X$ and the other to its negation.
We let $J$ be the set of these pairs.
By \autoref{NEQ-statement} there is an $([m],s)$-housing sequence of length at most $27$ and size at most $400s^{4/3}$ whose allocation set is $NEQ(J)$. Denote such a housing sequence by~$\Gcal_2$.

{\bf Definition of $\Gcal_3$.}
We let $K$ be the set of triplets of literals that form clauses in $F$.
By \autoref{OR-statement} there is an $([m],s)$-housing sequence of length at most $43$ and size at most $2001s^2$ whose allocation set is $OR(K)$. Denote such a housing sequence by $\Gcal_3$.

Having completed the definition of the $(I,m)$-housing sequences $\Gcal_1,\Gcal_2$ and $\Gcal_3$, we consider the hybrid deletion sequence $\Gcal=\Gcal_1\Gcal_2\Gcal_3$ formed by concatenation. 
It is an $([m],s)$-housing sequence of length at most $95$.
By \autoref{joint-allocation}, its allocation set $A$ is equal to:\vspace{-.25cm}
\[
A= EQ(I)\cap NEQ(J)\cap OR(K)\vspace{-.25cm}
\]

So $A$ contains those functions $f$ from $[m]$ to $\{0,1\}$ so that $f(i)=f(j)$ whenever $X_i=X_j$, and $f(i)\neq f(j)$ whenever $X_i=\neg X_j$, and $f(i)\vee f(j) \vee f(k)=1$ whenever $(X_i,X_j,X_k)$ forms a clause. Roughly speaking, $A$ contains those functions from $[m]$ to $\{0,1\}$ that \lq satisfy\rq\ $F$. Thus $F$ is satisfiable if and only if $A$ is nonempty, which is the case if and only if $\Gcal$ admits a simultaneous embedding. 
\end{pf}

\section{Basic constructions}\label{sec:part2}

In this section, we prove \autoref{EQ-statement}, \autoref{NEQ-statement} and \autoref{OR-statement} in a special case. Then in \autoref{sec:part3}, we deduce these lemmas from those special cases. Note that we often delimit deletion sequences by the (strict) subgraph relations $\subseteq$ and $\supseteq$ instead of commas to emphasise which graphs are minors of one another.

\begin{dfn}[Arches]
    Let $\Gcal$ be an $(I,s)$-housing sequence for some $s\in\mathbb{N}$ and $I\se [s]$. Given $i\in I$, we refer to the vertices $\{s_{4i-3},\ldots,s_{4i}\}$ in the surface as the planet as the \emph{foundation} of the house $H_i$. We say that $\Gcal$ a pair $(i,j)\in I^2$ is an \emph{arch} for $\Gcal$, if there is a path in some graph of~$\Gcal$ between the foundations of the houses $H_i$ and $H_j$ that avoids all other vertices of the planet.
\end{dfn}

\begin{figure}
    \begin{minipage}{0.54\textwidth}
        \large\scalebox{0.75}{\input{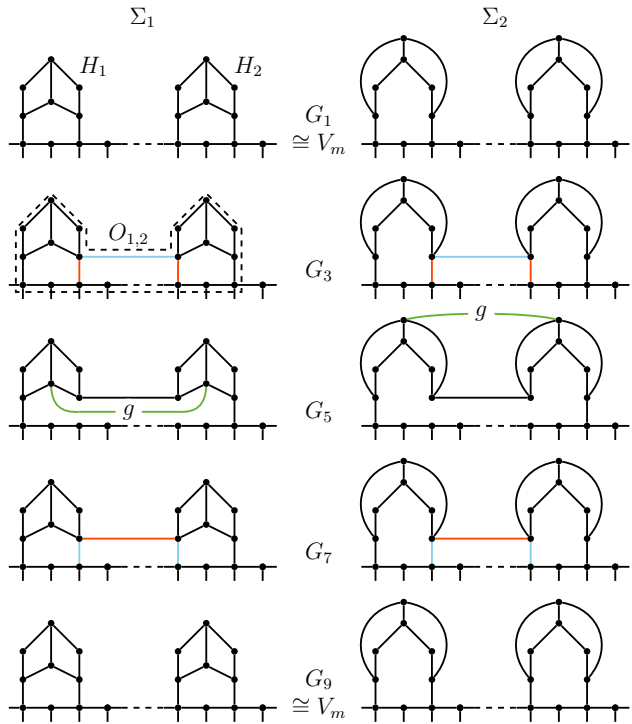}}
    \end{minipage}
    \hfill%
    \begin{minipage}{0.44\textwidth}
        \caption{The ``Equaliser gadget'' on houses: a housing sequence $(G_1\supseteq\ldots\subseteq G_9)$. Illustrated are the houses $H_1$ and $H_2$ in the graphs with odd index. The graphs of even index are implied: $G_2$ is equal to $G_1$, $G_4$ is obtained from $G_3$ by deleting the two \rcol\ edges, $G_6$ is obtained from $G_5$ by deleting the \gcol\ edge $g$ and $G_8$ is obtained from $G_9$ by deleting the orange edge. In general, \bcol\ edges are new edges, \rcol\ edges are edges which are to be deleted in the next step and \gcol\ edges are new edges which are to be deleted in the next step. Also labelled is the cycle $O_{1,2}$, which appears in all the graphs $G_3,\ldots,G_7$. The figure depicts the two (up to reorientation) simultaneous embeddings $\Sigma_1$ and $\Sigma_2$ of $(G_1\supseteq\ldots\subseteq G_9)$.}
    \label{fig:equaliser}
    \end{minipage}
\end{figure}

\begin{eg}
    The `Equaliser Gadget' in \autoref{fig:equaliser} has an arch $(1,2)$ since there is a path in $G_3$ from $v_1$ of $H_1$ to $v_2$ of $H_2$ (and also in $G_4,G_5,G_6,G_7$) using no other vertices of the planet.
\end{eg}

\begin{lem}\label{EQ-statement-special}
For every $m\in \Nbb$ and $(i,j)\in \binom{[m]}{2}$, there is an $(\{i,j\},m)$-housing sequence of length~$9$ and size at most $100\,m$ with allocation set $EQ(\{i,j\})$ and exactly one arch, namely $\{i,j\}$.
\end{lem}
\begin{pf}
    We will define a housing sequence $(G_1\supseteq\ldots\subseteq G_9)$ in order (illustrated in \autoref{fig:equaliser}):\vspace{-0.2cm}
    \begin{enumerate}
        \item  Firstly, let $G_1$ be the village of equator $4m$ and index set $\{i,j\}$.\vspace{-0.25cm}
        \item Let $G_2$ be a copy of $G_1$.\vspace{-0.25cm}
        \item Let $G_3$ be obtained from $G_2$ by adding an edge between $v_6$ of $H_i$ and $v_2$ of $H_j$, label it $g$.\vspace{-0.25cm}
        \item Let $G_4$ be obtained from $G_3$ by deleting the edge $v_6v_7$ from $H_i$ and the edge $v_1v_2$ from $H_j$.\vspace{-0.25cm}
        \item Let $G_5$ be obtained from $G_4$ by adding an edge between the vertices $w$ of the inhabitants of both houses.\vspace{-0.25cm}
        \item Finally, let $G_6,G_7,G_8$ and $G_9$ be copies of $G_4,G_3,G_2$ and $G_1$, respectively. In particular, $G_9$ is again the village of equator $4m$ and index set $\{i,j\}$.\vspace{-0.2cm}
    \end{enumerate}
    Let all of the edges in $(G_1\supseteq\ldots\subseteq G_9)$ be strict. We require the following claim.
    
    \begin{claim}\label{clm:in_in}
        The hybrid deletion sequence $(G_1\supseteq\ldots\subseteq G_9)$ admits exactly two simultaneous embeddings up to reorientation. In one of these, in both $G_1$ and $G_9$, both the inhabitants of $H_i$ and $H_j$ are embedded inside of their houses. In the other, both are embedded outside.
    \end{claim}
    \begin{cproof}
        Since all of the edges in $(G_1\supseteq\ldots\subseteq G_9)$ are strict, the simultaneous embeddings of $(G_1\supseteq\ldots\subseteq G_9)$ are exactly the simultaneous embeddings of the strict deletion sequence on the same graphs. Let $\Sigma$ be a simultaneous embedding of $(G_1\supseteq\ldots\subseteq G_9)$. 
        Let $O_{1,2}$ denote the cycle in $G_3,\ldots,G_7$ composed of the path $v_1\ldots v_6$ in $H_1$, the edge between $v_6$ in $H_1$ and $v_2$ in $H_2$, the path $v_2\ldots v_7$ in $H_2$ and the path in the surface of the planet from $v_1$ in $H_1$ to $v_7$ in $H_2$ in the positive direction. By the Jordan Curve Theorem, the cycle $O_{1,2}$ partitions $\Sbb^2$ into two regions in any embedding of one of the graphs $G_3,\ldots,G_7$. These regions form the two choices for where the edge $g$ and the inhabitants of the houses embed in an embedding of $G_5$. Each choice fixes a simultaneous embedding of $(G_1\supseteq\ldots\subseteq G_9)$. In one such simultaneous embedding, both inhabitants are embedded on the outside of their houses and in the other on the inside.
    \end{cproof}

    By the above claim, it follows that $(G_1\supseteq\ldots\subseteq G_9)$ is a housing sequence and it has allocation set $EQ(\{i,j\})$. By construction, it has one arch, namely $\{i,j\}$.
\end{pf}

\begin{figure}
    \begin{minipage}{0.54\textwidth}
        \large\scalebox{0.75}{\input{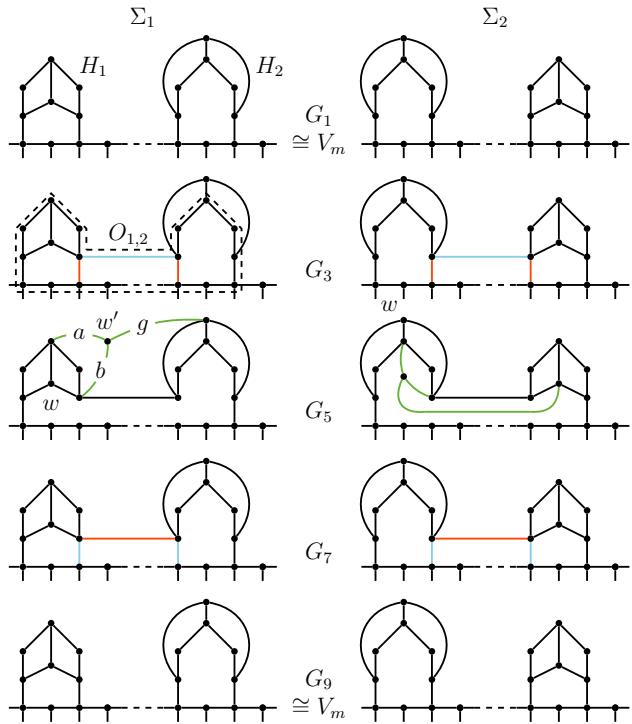}}
    \end{minipage}
    \hfill%
    \begin{minipage}{0.44\textwidth}
        \caption{The ``Negator gadget'' on houses: a housing sequence $(G_1\supseteq\ldots\subseteq G_9)$. Illustrated are the houses $H_1$ and $H_2$ in the graphs with odd index. The graphs of even index are implied: $G_2$ is equal to $G_1$, $G_4$ is obtained from $G_3$ by deleting the two \rcol\ edges, $G_6$ is obtained from $G_5$ by deleting the vertex $w'$ and the \gcol\ edges $a,b,g$ and $G_8$ is obtained from $G_9$ by deleting the \rcol\ edge. In general, \bcol\ edges are new edges, \rcol\ edges are edges which are to be deleted in the next step and \gcol\ edges are new edges which are to be deleted in the next step. Also labelled is the cycle $O_{1,2}$, which appears in all the graphs $G_3,\ldots,G_7$. The figure depicts the two (up to reorientation) simultaneous embeddings $\Sigma_1$ and $\Sigma_2$ of $(G_1\supseteq\ldots\subseteq G_9)$.}
    \label{fig:negator}
    \end{minipage}
\end{figure}

\begin{lem}\label{NEQ-statement-special}
For every $m\in \Nbb$ and $(i,j)\in \binom{[m]}{2}$, there is an $(\{i,j\},m)$-housing sequence of length 9 and size at most $100\,m$ with allocation set $NEQ(\{i,j\})$ and exactly one arch, namely~$\{i,j\}$.
\end{lem}
\begin{pf}
    We will define a housing sequence $(G_1\supseteq\ldots\subseteq G_9)$ in order (illustrated in \autoref{fig:negator}):\vspace{-0.2cm}
    \begin{enumerate}
        \item  Firstly, let $G_1$ be the village of equator $4m$ and index set $\{i,j\}$.\vspace{-0.25cm}
        \item Let $G_2$ be a copy of $G_1$.\vspace{-0.25cm}
        \item Let $G_3$ be obtained from $G_2$ by adding an edge between $v_6$ of $H_i$ and $v_2$ of $H_j$.\vspace{-0.25cm}
        \item Let $G_4$ be obtained from $G_3$ by deleting the edge $v_6v_7$ from $H_i$ and the edge $v_1v_2$ from $H_j$.\vspace{-0.25cm}
        \item Let $G_5$ be obtained from $G_4$ by gluing the leaves of a new claw to $v_4$ and $v_6$ in $H_i$ and the inhabitant $w$ in $H_j$. Call these new edges $a$,$b$ and $g$, respectively. We will denote the new vertex by $w'$ and call it the \emph{exhabitant} of $H_i$.\vspace{-0.25cm}
        \item Finally, let $G_6,G_7,G_8$ and $G_9$ be copies of $G_4,G_3,G_2$ and $G_1$, respectively. In particular, $G_9$ is again the village of equator $4m$ and index set $\{i,j\}$.\vspace{-0.2cm}
    \end{enumerate}
    Let all of the edges in $(G_1\supseteq\ldots\subseteq G_9)$ be strict.
    
    Let $O_{1,2}$ denote the cycle in $G_3,\ldots,G_7$ composed of the path $v_1\ldots v_6$ in $H_1$, the edge between $v_6$ in $H_1$ and $v_2$ in $H_2$, the path $v_2\ldots v_7$ in $H_2$ and the path in the surface of the planet from $v_1$ in $H_1$ to $v_7$ in $H_2$ in the positive direction. By the Jordan Curve Theorem, the cycle $O_{1,2}$ partitions $\Sbb^2$ into two regions in any embedding of one of the graphs $G_3,\ldots,G_7$. We call the region containing the centre of the planet the \emph{outside} of $O_{1,2}$ and the other region we call the \emph{inside}.
    
    \begin{claim}\label{clm:in_ex}
        In every embedding of $G_5$ the inhabitant and exhabitant of $H_i$ are embedded in different regions bounded by the cycle $O_{1,2}$.
    \end{claim}
    \begin{cproof}
        Let $H$ be the minor of $G_5$ obtained by contracting all edges except those spanning the pairs $v_1v_2$, $v_2v_3$ and $v_4v_5$ in $H_i$ and the edges incident with either the inhabitant $w$ or exhabitant $w'$ of $H_i$. The graph $H$ consists of a $3$-cycle (which is a minor of $O_{1,2}$) along with two claws whose centres are $w$ and $w'$. Observe that if $w$ and $w'$ are embedded in the same region bounded by $O_{1,2}$ in an embedding $\iota$ of $G_5$, then $w$ and $w'$ are both embedded in the same side of the $3$-cycle in the embedding induced by $\iota$ on $H$.  However, such an embedding is impossible.
    \end{cproof}
    \begin{claim}
        The hybrid deletion sequence $(G_1\supseteq\ldots\subseteq G_9)$ admits exactly two simultaneous embeddings up to reorientation. In one of these, in both $G_1$ and $G_9$ the inhabitant of $H_i$ is embedded on the inside of its house and the inhabitant of $H_j$ on the outside of its house. In the other, the roles of `$i$' and '$j$' are swapped.
    \end{claim}
    \begin{cproof}
         Since all of the edges in $(G_1\supseteq\ldots\subseteq G_9)$ are strict, the simultaneous embeddings of $(G_1\supseteq\ldots\subseteq G_9)$ are exactly the simultaneous embeddings of the strict deletion sequence on the same graphs. Let $\Sigma$ be a simultaneous embedding of $(G_1\supseteq\ldots\subseteq G_9)$. The inside and outside of $O_{1,2}$ form two possibilities for where the edge $g$ and its endvertices ($w'$ in $H_i$ and $w$ in $H_j$) embed. Each of these choices fixes a simultaneous embedding of $(G_1\supseteq\ldots\subseteq G_9)$.
         
         By \autoref{clm:in_ex}, we get that the exhabitant of $H_i$ and the inhabitant of $H_j$ must lie in distinct regions bounded by the cycle $O_{1,2}$. This proves our claim.
    \end{cproof}

    By the above claim, it follows that $(G_1\supseteq\ldots\subseteq G_9)$ is a housing sequence and has allocation set $NEQ(\{i,j\})$. By construction, it has one arch, namely $\{i,j\}$.
\end{pf}

In the proofs of \autoref{EQ-statement-special} and \autoref{NEQ-statement-special}, we force all of the edges in our deletion sequences to be strict, so that they behave exactly as strict deletion sequences. In fact, we confine all of the `weak' behaviour used for the proof of \autoref{thm:main} to proof of the following \autoref{OR-statement-special}.

\begin{figure}
    \begin{minipage}{0.50\textwidth}
        \large\scalebox{0.9}{\input{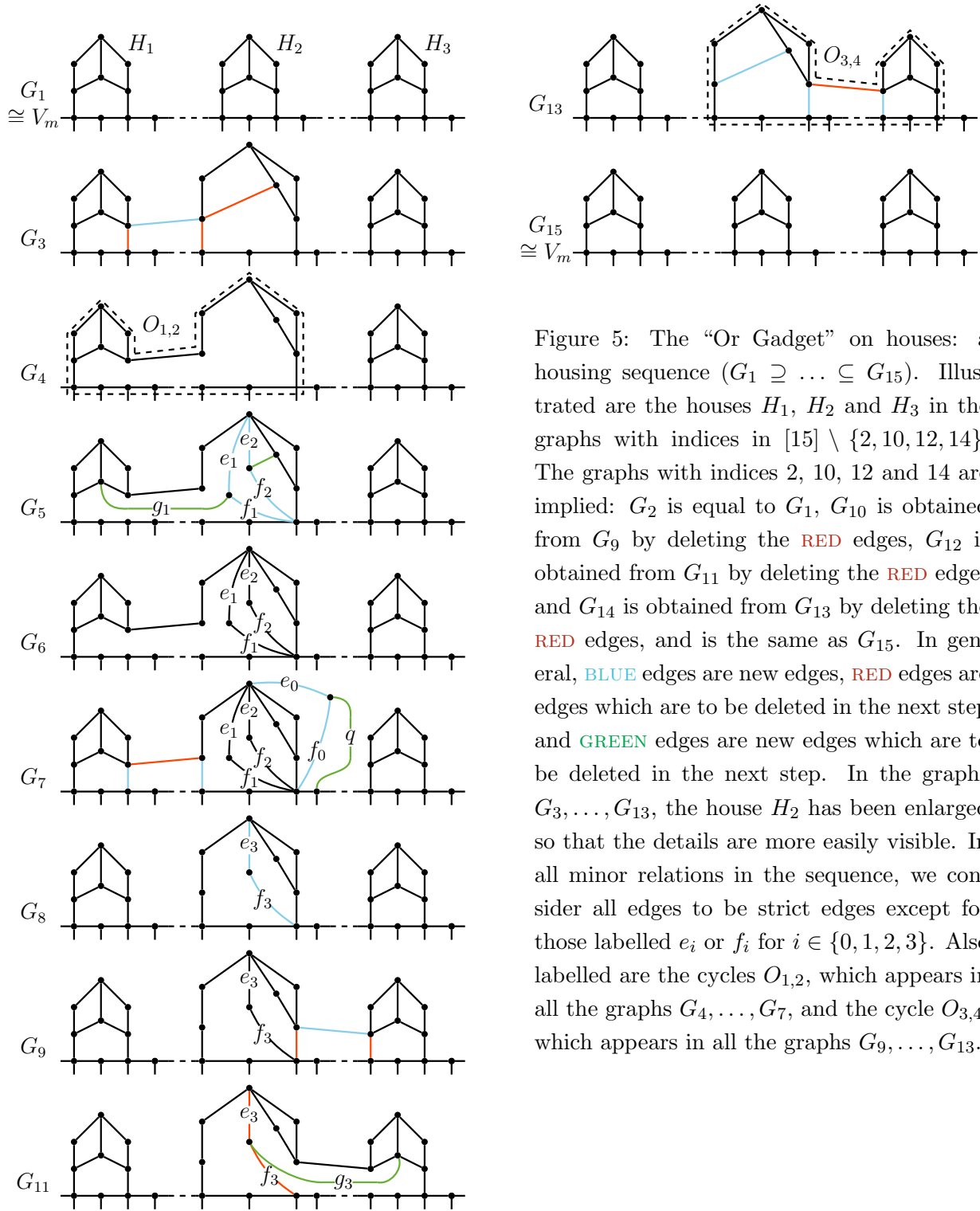}}
    \end{minipage}
    \hfill%
    \begin{minipage}{0.46\textwidth}
        \vspace{2.5cm}
        \caption{The ``Or Gadget'' on houses: a housing sequence $(G_1\supseteq\ldots\subseteq G_{15})$. Illustrated are the houses $H_1$, $H_2$ and $H_3$ in the graphs with indices in $[15]\sm\{2,10,12,14\}$. The graphs with indices $2$, $10$, $12$ and $14$ are implied: $G_2$ is equal to $G_1$, $G_{10}$ is obtained from $G_9$ by deleting the \rcol\ edges, $G_{12}$ is obtained from $G_{11}$ by deleting the \rcol\ edges and $G_{14}$ is obtained from $G_{13}$ by deleting the \rcol\ edges, and is the same as $G_{15}$. In general, \bcol\ edges are new edges, \rcol\ edges are edges which are to be deleted in the next step and \gcol\ edges are new edges which are to be deleted in the next step. In the graphs $G_3,\ldots,G_{13}$, the house $H_2$ has been enlarged so that the details are more easily visible. In all minor relations in the sequence, we consider all edges to be strict edges except for those labelled $e_i$ or $f_i$ for $i\in\{0,1,2,3\}$. Also labelled are the cycles $O_{1,2}$, which appears in all the graphs $G_4,\ldots,G_7$, and the cycle $O_{3,4}$, which appears in all the graphs $G_9,\ldots,G_{13}$.}
    \label{fig:fror}
    \end{minipage}
\end{figure}

\begin{lem}\label{OR-statement-special}
For every $m\in 3\Nbb$ and $i\in [m/3]$, there is an $(3i-2,3i-1,3i)$-housing sequence of length $15$ and size at most $300\,m$ with allocation set $OR(\{3i-2,3i-1,3i\})$ and whose set of arches is $\{(3i-2,3i-1),(3i-1,3i)\}$.
\end{lem}
\begin{pf}
    We will define a housing sequence $(G_1\supseteq\ldots\subseteq G_{15})$ in order (illustrated in \autoref{fig:fror}):\vspace{-0.2cm}
    \begin{enumerate}
        \item Firstly, let $G_1$ be the village with equator $m$ and index set $\{1,2,3\}$.\vspace{-0.25cm}
        \item Let $G_2$ be a copy of $G_1$.\vspace{-0.25cm}
        \item Let $G_3$ be obtained from $G_2$ by adding an edge between $v_6$ of $H_1$ and $v_2$ of $H_2$.\vspace{-0.25cm}
        \item Let $G_4$ be obtained from $G_3$ by deleting the edge spanning the pair $v_6v_7$ in $H_1$ and the edges spanning the pairs $v_1v_2$ and $v_2w$ in $H_2$.\vspace{-0.25cm}
        \item Let $G_5$ be obtained from $G_4$ by adding the following structure: Let $Q_1$ and $Q_2$ denote paths of length $2$ on edges labeled $e_1f_1$ and $e_2f_2$, respectively. Identify the starts of $Q_1$ and $Q_2$ with $v_4$ in $H_2$ and identify the ends of $Q_1$ and $Q_2$ with $v_7$ in $H_2$. Connect the internal vertex of the path $Q_1$ to the inhabitant of $H_1$ via a new edge $g_1$ and connect the internal vertex of the path $Q_2$ to the inhabitant of $H_2$ via a new edge $g_2$.\vspace{-0.25cm}
        \item Let $G_6$ be obtained from $G_5$ by deleting $g_1$ and $g_2$.\vspace{-0.25cm}
        \item Let $G_7$ be obtained from $G_6$ by adding edges spanning the pairs $v_6v_7$ in $H_1$ and $v_1v_2$ in $H_2$ along with the following structure: Let $Q_0$ be a path with edges labelled $e_0f_0$. Identify the start of $Q_0$ with $v_4$ in $H_2$ and the end of $Q_0$ with $v_7$ in $H_2$. Then connect the interior vertex of $Q_0$ to the first vertex along from the house $H_2$ on the surface of the planet in the positive direction. Call this new edge $q$.\vspace{-0.25cm}
        \item Let $G_8$ be the graph obtained from $G_7$ by deleting the edges $\{e_0,f_0,\;e_1,f_1,\;e_2,f_2,\;q\}$ and adding the following structure: Let $Q_3$ be the path of length $2$ with edges labelled $e_3f_3$ and identify the start of $Q_3$ with the vertex $v_4$ in $H_2$ and the end of $Q_3$ with the vertex $v_7$ in~$H_2$.\vspace{-0.25cm}
        \item Let $G_9$ be the graph obtained from $G_8$ by adding from $v_6$ in $H_2$ to $v_2$ in $H_3$.\vspace{-0.25cm}
        \item Let $G_{10}$ be the graph obtained from $G_9$ be deleting the edge spanning $v_6v_7$ in $H_2$ and the edge spanning $v_1v_2$ in $H_3$..\vspace{-0.25cm}
        \item Let $G_{11}$ be the graph obtained from $G_{10}$ by adding an edge labelled $g_3$ between the interior vertex of $Q_3$ and the inhabitant of $H_3$.\vspace{-0.25cm}
        \item Let $G_{12}$ be the graph obtained from $G_{11}$ by deleting the edges $e_3,f_3$ and $g_3$ along with the interior vertex of $Q_3$.\vspace{-0.25cm}
        \item Let $G_{13}$ be the graph obtained from $G_{12}$ by adding edges spanning the pairs $v_1w$ and $v_6v_7$ in $H_2$ and the pair $v_1v_2$ in $H_3$.\vspace{-0.25cm}
        \item Let $G_{14}$ be obtained from $G_{15}$ by deleting the edge between $v_6$ in $H_2$ and $v_2$ in $H_3$.\vspace{-0.25cm}
        \item Let $G_{15}$ be a copy of $G_{14}$ and note that $G_{15}$ is again the village with equator $m$ and index set $\{1,2,3\}$.\vspace{-0.25cm}
    \end{enumerate}
    Let all of the edges in $(G_1\supseteq\ldots\subseteq G_{15})$ be strict except for those labelled $e_i$ or $f_i$ for $i\in\{0,1,2,3\}$. Let $A$ denote the strict edges.

    Let $O_{1,2}$ denote the cycle in $G_3,\ldots,G_7$ composed of the path $v_1\ldots v_6$ in $H_1$, the edge between $v_6$ in $H_1$ and $v_2$ in $H_2$, the path $v_2\ldots v_7$ in $H_2$ and the path in the surface of the planet from $v_1$ in $H_1$ to $v_7$ in $H_2$ in the positive direction. By the Jordan Curve Theorem, the cycle $O_{1,2}$ partitions $\Sbb^2$ into two regions in any embedding of one of the graphs $G_3,\ldots,G_7$. We call the region containing the centre of the planet the \emph{outside} of $O_{1,2}$ and the other region we call the \emph{inside}.

    We make the analogous definitions for the cycle $O_{2,3}$ in the graphs $G_9,\ldots,G_{13}$ by writing `$H_2$' in the place of `$H_1$' and `$H_3$' in the place of `$H_2$'. We require the following claims.

    \begin{claim}\label{clm:housing}
        The hybrid deletion sequence $(G_1\supseteq\ldots\subseteq G_{15})$ is a housing sequence.
    \end{claim}
    \begin{cproof}
        In any simultaneous embedding, observe that for each $i\in\{1,2,3\}$, the inhabitant of $H_i$ in $G_1$ is embedded inside its house if and only if it is embedded inside its house in $G_{15}$. In particular, we note that the inhabitant of $H_1$ (resp. $H_2$) is embedded inside of its house in $G_1,G_2,G_3$ if and only if it is embedded inside of the cycle $O_{1,2}$ in $G_4,\ldots,G_7$ if and only if it is embedded inside of its house in $G_8$. The inhabitant of $H_3$ is embedded inside of its house in $G_1$ if and only if it is embedded inside of its house in $G_8$. For the graphs $G_8,\ldots,G_{15}$, similar reasoning applies.
    \end{cproof}

    \begin{claim}\label{clm:control12}
        Let $i\in\{1,2\}$. In any simultaneous embedding of $(G_1\supseteq\ldots\subseteq G_{15})$, the path $Q_i$ is embedded inside of $O_{1,2}$ in $G_5$ if and only if the inhabitant of $H_i$ is embedded inside its house in~$G_1$.
    \end{claim}
    \begin{cproof}
        We observe that the inside of the $O_{1,2}$ in $G_3$ contains the inside of both of the houses $H_1$ and $H_2$.
    \end{cproof}

    \begin{claim}\label{clm:q_one}
        Let $\Sigma$ be a simultaneous embedding of $(G_1\supseteq\ldots\subseteq G_{15})$. Let $q$ denote the number of paths from $\{Q_0,Q_1,Q_2\}$ embedded inside of $O_{1,2}$ in $G_7$ under $\Sigma$. Then $q$ is equal to the number of inhabitants from the houses $H_1$ and $H_2$ embedded inside of their houses in $G_1$ under $\Sigma$.
    \end{claim}
    \begin{cproof}
        Let $\sigma_5,\sigma_6$ and $\sigma_7$ be embeddings of $G_5$, $G_6$ and $G_7$, respectively, and let $q_5$ be the number of paths from $\{Q_1,Q_2\}$ that are embedded on the inside of $O_{1,2}$ in $\sigma_5$. Define $q_6$ and $q_7$ analogously but with `$5$' replaced by `$6$' and `$5$' replaced by `$7$', respectively. Suppose that, up to $A$-fixing isomorphism, we have that $\sigma_6$ is the embedding induced on $G_6$ by $G_5$. In this case, the paths $\{Q_1,Q_2\}$ in $G_5$ are mapped to themselves in $G_6$ but perhaps not via the identity. However, either possible map yields that $q_6=q_5$. Assuming that $\sigma_6$ is, up to $A$-fixing isomorphism, an embedding induced on $G_6$ by $\sigma_7$, we similarly get that $q_7=q_6$. Now applying \autoref{clm:control12} completes the proof of our claim.
    \end{cproof}

    \begin{claim}\label{clm:control34}
        In any simultaneous embedding of $(G_1\supseteq\ldots\subseteq G_{15})$, the path $Q_3$ is embedded inside of $O_{3,4}$ in $G_9$ if and only if the inhabitant of $H_3$ is embedded inside its house in $G_{15}$.
    \end{claim}
    \begin{cproof}
        analogous to the proof of \autoref{clm:control12}.
    \end{cproof}

    \begin{claim}\label{clm:controlh2}
        In any simultaneous embedding of $(G_1\supseteq\ldots\subseteq G_{15})$, the path $Q_3$ is embedded inside the house $H_2$ in $G_8$ if and only if the inhabitant of $H_3$ is embedded inside its house in $G_{15}$.
    \end{claim}
    \begin{cproof}
        This follows immediately from \autoref{clm:control34}.
    \end{cproof}

    \begin{claim}\label{clm:q_null}
        Let $\Sigma$ be a simultaneous embedding of $(G_1\supseteq\ldots\subseteq G_{15})$. Let $q$ denote the number of paths from $\{Q_0,Q_1,Q_2\}$ embedded inside of $O_{1,2}$ in $G_7$ under $\Sigma$. If $q=0$, then the inhabitant of $H_3$ is embedded on the outside of its house in $G_{15}$ under $\Sigma$.
    \end{claim}
    \begin{cproof}
        Suppose that $q=0$. Let $\sigma_7$ and $\sigma_8$ denote the embeddings that $\Sigma$ assigns to $G_7$ and $G_8$, respectively. Then $\sigma_8$ is, up to $A$-fixing isomorphism, an embedding induced on $G_8$ by $\sigma_7$. Now observe that our $A$-fixing isomorphism must sent one of the paths $Q_0,Q_1,Q_2$ to $Q_3$. Hence $Q_3$ cannot appear inside of the house $H_2$ in $G_8$. Applying \autoref{clm:control34} completes the proof of our claim.
    \end{cproof}

    Combining \autoref{clm:q_null} with \autoref{clm:housing} proves that there is no simultaneous embedding of $(G_1\supseteq\ldots\subseteq G_{15})$ so that all three inhabitants of the houses $H_1,H_2,H_3$ are embedded outside of their houses in $G_1$. It remains to show that every other combination of being inside/outside of houses is represented by some simultaneous embedding. However, one may check this easily by hand. This completes the proof of the lemma.
\end{pf}

\section{Proofs of the lemmas \autoref{EQ-statement}, \autoref{NEQ-statement} and \autoref{OR-statement}}\label{sec:part3}

In this section, we deduce lemmas \autoref{EQ-statement}, \autoref{NEQ-statement} and \autoref{OR-statement} from the special cases proved in the last subsection. In the first subsection, we explain overarching ideas and then perform the individual proofs in the next subsection.

\subsection{Overarching ideas}
\label{sec:overarching}
Let $\Gcal=(G_1,\ldots,G_n)$ and $\Hcal=(H_1,\ldots,H_n)$ be two hybrid deletion sequences of the same length $n$. We refer to the sequence $(G_1\cup H_1,\ldots,G_n\cup H_n)$, where $G_i\cup H_i$ is the graph $(V(G_i)\cup V(H_i), E(G_i)\cup E(H_i))$, as the \emph{union} of $\Gcal$ and $\Hcal$. We stress that vertex labels matter for this definition, and there can be vertices in $V(G_i)\cap V(H_i)$. We denote the union by $\Gcal\cup \Hcal$.

\begin{lem}\label{union-lem}
Let $\Gcal$ and $\Hcal$ be housing sequences of the same length $n$ with the same equator $m$ and disjoint index sets $I$ and $J$. Then $\Gcal\cup \Hcal$ is a housing sequence of length $n$, equator $m$ and index set $I\cup J$.\qed
\end{lem}

We say that two graph $G$ and $H$ are \emph{unifiable} if for every embedding $\sigma$ of $G$ and every embedding $\pi$ of $H$ there is a unique embedding of $G\cup H$ that induces $\sigma$ and $\pi$.
In these circumstances we shall denote this unique embedding of $G\cup H$ by $\sigma\cup \pi$. 
We say that two hybrid deletion sequences $\Gcal=(G_1,\ldots,G_n)$ and $\Hcal=(H_1,\ldots,H_n)$ of the same length $n$ are \emph{unifiable} if for all $i\in [n]$ we have that $G_i$ and $H_i$ are unifiable, and furthermore the set of weak edges of $\Gcal$ is disjoint from the set of weak edges of $\Hcal$. 
Let $\Sigma=(\sigma_1,\ldots,\sigma_n)$ and $\Pi=(\pi_1,\ldots,\pi_n)$ be simultaneous embeddings of unifiable hybrid deletion sequences $\Gcal$ and $\Hcal$, their union is $\Sigma\cup \Pi= (\sigma_1\cup\pi_1,\ldots,\sigma_n\cup \pi_n)$. 

\begin{lem}\label{unification-lem}
The union $\Sigma\cup \Pi$ of two simultaneous embeddings $\Sigma$ and $\Pi$ of unifiable hybrid deletion sequences $\Gcal$ and $\Hcal$, respectively, is a simultaneous embedding of $\Gcal\cup \Hcal$.
\end{lem}
\begin{pf}
Let $H$ and $G$ be graphs and suppose that $H$ is a hybrid minor of $G$ with respect to $A\subseteq E(G)$, the set of strict edges. We say that an embedding $\iota$ of $G$ \emph{weakly induces} an embedding $\iota'$ of $H$ if $\iota'$ is the same up to $A$-fixing isomorphism as the embedding induced by $\iota$ on some subgraph of $G$. Let $i$ be an odd integer. We are to show that the embedding $\sigma_i\cup \pi_i$ weakly induces the embedding $\sigma_{i-1}\cup \pi_{i-1}$ (whenever $i>1$) and the embedding $\sigma_{i+1}\cup \pi_{i+1}$ (whenever $i<n$). 
So assume that $i>1$ and we shall show that $\sigma_i\cup \pi_i$ weakly induces $\sigma_{i-1}\cup \pi_{i-1}$.

Let $S$ and $P$ denote the sets of weak edges in the hybrid minor relations between $H_{i-1}$ and $H_i$ and between $G_{i-1}$ and $G_i$, respectively. Since $\Sigma$ and $\Pi$ are hybrid simultaneous embeddings, the embeddings $\sigma_i$ and $\pi_i$ weakly induce the embeddings $\sigma_{i-1}$ and $\pi_{i-1}$ on $G_{i-1}$ and $H_{i-1}$, respectively. So there exist subgraphs $G_{i-1}'$ and $H_{i-1}'$ of $G_i$ and $H_i$, respectively, and isomorphisms $g:G_{i-1}'\rightarrow G_{i-1}$ and $h:H_{i-1}'\rightarrow H_{i-1}$ which are $E(G_i)\sm S$-fixing and $E(H_i)\sm P$-fixing, respectively, witnessing this fact. By assumption, $S\cap P=\varnothing$. Hence the graph $H_{i-1}'\cup G_{i-1}'$ is well-defined. Let $k:G_{i-1}'\cup H_{i-1}'\rightarrow G_{i-1}\cup H_{i-1}$ be the unique isomorphism restricting to both $g$ and $h$. Let $\phi_{i-1}$ be the image under $k$ of the embedding induced on $G_{i-1}'\cup H_{i-1}'$ by $\sigma_i\cup\pi_i$. Observe that $\phi_{i-1}$ induces both $\sigma_{i-1}$ on $G_{i-1}$ and $\pi_{i-1}$ on $H_{i-1}$. Since $G_{i-1}$ and $H_{i-1}$ are unifiable, the uniqueness condition implies that $\phi_{i-1}=\sigma_{i-1}\cup\pi_{i-1}$.
Hence $k$ witnesses that $\sigma_i\cup \pi_i$ weakly induces $\sigma_{i-1}\cup \pi_{i-1}$.
One proves similarly that $\sigma_i\cup \pi_i$ weakly induces $\sigma_{i+1}\cup \pi_{i+1}$ (whenever $i<n$). 
\end{pf}

\begin{eg}
Given two housing sequences $\Gcal$ and $\Hcal$ of same length and equator, the set of arches of $\Gcal\cup \Hcal$ is the union of the sets of arches of $\Gcal$ and $\Hcal$. 
\end{eg}

We say that two pairs $\{a,b\},\{c,d\}\in \binom{m}{2}$ cross if when we order $a$, $b$, $c$, $d$ by cyclic order of the equator $\Zbb_m$, then this sequence alternates between the pairs $\{a,b\}$ and $\{c,d\}$. We say that two arches \emph{cross} when their underlying pairs cross.

Our next goal is to prove the following lemma.

\begin{lem}\label{are-unifiable}
Let $\Gcal$ and $\Hcal$ be two housing sequences of the same length and the same equator with disjoint index sets and disjoint sets of weak edges. If no arch from $\Gcal$ crosses an arch from $\Hcal$, then $\Gcal$ and $\Hcal$ are unifiable.
\end{lem}

We shall prepare for an inductive proof with the following being the base case.

We say that a pair $(a,b)$ of elements of an equator of a housing sequence \emph{separates} a set $J$ of indices if there are $c,d$ vertices in the foundations of the houses $J$ so that the sequence $a$, $b$, $c$, $d$ alternates in the cyclic ordering of the equator. And we say that two index sets \emph{separate each other} if a pair of foundation vertices from one separates the other. 

\begin{lem}\label{are-unifiable-base}
Let $\Gcal$ and $\Hcal$ be two housing sequences of the same length and the same equator with disjoint index sets $I$ and $J$, respectively, and disjoint sets of weak edges. If $I$ and $J$ do not separate each other, then $\Gcal$ and $\Hcal$ are unifiable.\qed
\end{lem}
\begin{pf}
Let $G_k$ and $H_k$ denote the $k$-th graphs in $\Gcal$ and $\Hcal$, respectively. Since $I$ and $J$ do not separate each other, there exists a pair $(a,b)$ of equator vertices whose removal splits the equator into two arcs, each containing only foundation vertices of the houses for $I$ and $J$, respectively. Let $P$ be the path of length $2$ joining $a$ and $b$ that uses the center of the planet and observe that removing $P$ separates both $G_k$ and $H_k$ into two components. Denote by $G_k'$ and $H_k'$ the subgraphs of $G_k$ and $H_k$, respectively, obtained by deleting everything on the side of $P$ not containing the foundations of the houses for $J$ and $I$, respectively. Embedding $G_k'$ and $H_k'$ into two corresponding half-planes, with $P$ in the boundary between them, yields an embedding of $G_k\cup H_k$. This proves that they are unifiable.
\end{pf}

Given a housing sequence with index set $I$, a \emph{break} is a partition of $I$ into two nonempty subsets such that there is no arch with an endpoint in each of these two sets, and these sets do not separate one another.

\begin{lem}\label{break}
Let $G$ be a housing sequence with equator $m$ and index set $I$, and suppose $I$ admits a break $I=I_1\dot\cup I_2$. Then $G$ can be written as a union of two housing sequences with the same equator $m$ and index sets $I_1$ and $I_2$, respectively.
\end{lem}
\begin{pf}
For $r\in\{1,2\}$, let $G{\upharpoonright I_r}$ be the housing sequence obtained from $G$ by, at each step, deleting all vertices that are not joined to a foundation for $I_r$ by a path which internally avoids the planet. Because $I=I_1\dot\cup I_2$ is a break, no arch has one endpoint in each part, so the arches of $G$ split as a disjoint union of the arches of $G\upharpoonright I_1$ and $G\upharpoonright I_2$. Moreover, the break condition includes that $I_1$ and $I_2$ do not separate one another, so by \autoref{are-unifiable-base} the housing sequences $G\upharpoonright I_1$ and $G\upharpoonright I_2$ are unifiable. Hence, by \autoref{union-lem}, their union is a housing sequence on equator $m$ with index set $I_1\cup I_2=I$, and since the arch sets partition, this union is exactly $G$.
\end{pf}

\begin{lem}\label{create_a_break}
Let $\Gcal$ and $\Hcal$ be two housing sequences of the same length and the same equator with disjoint index sets $I$ and $J$, respectively. If there is an arch $(a,b)$ with both endpoints in $I$ that separates $J$, then $J$ admits a break, or $(a,b)$ crosses an arch from $\Hcal$.
\end{lem}
\begin{pf}
The indices $a$ and $b$ when removed from the equator leave two intervals. This defines a partition of $J$ into those indices in the first of these intervals versus the partition class consisting of those indices in the second of these intervals. If there is an arch between these two partition classes of $J$---so an arch with one endpoint in each of them---then $(a,b)$ crosses an arch from $\Hcal$. Otherwise this partition of $J$ is a break.
\end{pf}

\begin{pf}[Proof of \autoref{are-unifiable}]
Denote the index sets of $\Gcal$ and $\Hcal$ by $I$ and $J$, respectively. 
We colour the indices on the equator that are in foundations for houses in $I$ \rcol, and we colour those in $J$ \bcol.
We prove the statement by induction on the number of colour changes. If there are at most two colour changes, then $I$ and $J$ do not separate one another and \autoref{are-unifiable-base} completes the proof.

So assume that are at least four colour changes (note the number of changes is always an even number). Pick a \rcol\ index and take a maximal subinterval of the equator which contains it and contains no \bcol\ indices, denote it by $I'$. Either $(I',I\sm I')$ is a break or there is an arch between an index in $I'$ and an index in $I\sm I'$. By maximality of $I'$, this arch separates $J$. Since no arch of $I$ crosses an arch of $J$ by assumption, by \autoref{create_a_break}, $J$ admits a break. To summarise, $I$ or $J$ admits a break.

By symmetry assume that $I$ admits a break. Denote that break by $(I_1,I_2)$.
By \autoref{break} write $\Gcal$ as the union of two housing sequences with the same equator with index sets $I_1$ and $I_2$. Now we consider $\Hcal$ and $\Gcal_1$. This pair has disjoint index sets whose arches do not cross one another, and it has less colour changes on the equator than the pair $\Hcal$ and $\Gcal$. So $\Hcal$ and $\Gcal_1$ are unifiable by induction. Furthermore, the pair $\Hcal\cup \Gcal_1$ and $\Gcal_2$ has disjoint index sets whose arches do not cross one another, and it has less colour changes on the equator than the pair $\Hcal$ and $\Gcal$. Thus the pair $\Hcal\cup \Gcal_1$ and $\Gcal_2$ is unifiable and furthermore its union is equal to $\Gcal\cup \Hcal$, proving that $\Gcal$ and $\Hcal$ are unifiable.
\end{pf}



\begin{rem}
The converse of \autoref{are-unifiable} is also true but we do not need it in our proofs, so do not include a proof. 
\end{rem}

The \lq overarching ideas\rq\ of this subsection are applied in the forthcoming subsection to reduce the general situation step by step to the special cases from \autoref{sec:part2}. Structurally, each proof proceeds by decomposing along breaks or unions so that the arches involved do not cross.

\subsection{The proofs of the three lemmas}
\label{subsec:two}
First we shall derive the following special case of \autoref{EQ-statement} from \autoref{EQ-statement-special}.

\begin{lem}\label{EQ-statement-2}
For every $m\in \Nbb$, $s\geq m$ and set $P\se \binom{[m]}{2}$ of disjoint pairs so that no two of them cross, there is an $([m],s)$-housing sequence of length at most $9$ and size at most $100\,m^2$ with allocation set~$EQ(P)$. 
\end{lem}
\begin{pf}
    We will prove this by induction on the number of pairs in $P$. If $P$ consists of one pair, we simply apply \autoref{EQ-statement-special} and are done. For the inductive hypothesis, suppose that we have a $([m],s)$-housing sequence $\Gcal$ of length $9$ and size at most $100\,|P|\, m$ with allocation set $EQ(P)$. Further suppose that the set of arches of $\Gcal$ is exactly $P$. Let $P'=P\sqcup\{(i,j)\}$ for some additional pair $(i,j)\in \binom{[m]}{2}\sm P$. We would like to modify $\Gcal$ so that the allocation set becomes $EQ(P')$. To this end, let $\Gcal_{i,j}$ denote the $(\{i,j\},m)$-housing sequence given by \autoref{EQ-statement-special}; it has length 9 and size at most $100\,m$ with allocation set $EQ(\{i,j\})$ and exactly one arch, namely $\{i,j\}$. Further observe that the set of weak edges of $\Gcal$ is disjoint from the set of weak edges of $\Gcal_{i,j}$ (in fact, neither of them have weak edges).
    By \autoref{are-unifiable}, we have that $\Gcal$ and $\Gcal_{i,j}$ are unifiable, and denote their unification by $\Gcal':=\Gcal\cup \Gcal_{i,j}$.

    Now observe that $\Gcal'$ has size at most $100\,(|P|+1)\, m$ and length $9$ and its set of arches is exactly~$P'$. It remains to show that the allocation set for $\Gcal'$ is $EQ(P')$. Given a hybrid simultaneous embedding of $\Gcal'$, this induces a hybrid simultaneous embedding on both $\Gcal$ and $\Gcal_{i,j}$, and hence satisfies both $EQ(P)$ and $EQ(\{i,j\})$, thus satisfying $EQ(P')$. Conversely, given an arbitrary allocation the allocation function $f$ from $EQ(P')$, we have both $f\in EQ(P)$ and $f\in EQ(\{i,j\})$. Hence there are hybrid simultaneous embeddings $\Sigma$ and $\Sigma_{i,j}$ of $\Gcal$ and $\Gcal_{i,j}$, respectively, both satisfying $f$. We apply \autoref{unification-lem} to obtain a simultaneous embedding $\Sigma\cup\Sigma_{i,j}$ of $\Gcal'$ which satisfies~$f$. Hence the allocation set of $\Gcal'$ is $EQ(P')$.
\end{pf}

\begin{figure}
    
        {\large\scalebox{0.8}{\input{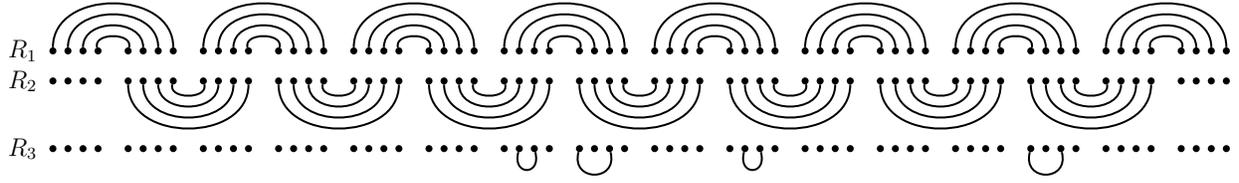}}}

        \caption{A depiction of sets of pairs $R_1,R_2$ and $R_3$, defined in the proof of \autoref{EQ-statement}, which are subsets of $\binom{[4^3]}{2}$. Here $P=\{\{2,3\},\{2,4\}\}$ is a subset of $\binom{[4]}{2}$. The transitive closure of the pairs from $R_1$, $R_2$ and $R_3$ amounts to `pairing up' $2$ with $3$ and $2$ with $4$ (when we restrict our attention to the first $4$ of $4^3$ points).}
    \label{fig:arching}
\end{figure}

\begin{pf}[Proof of \autoref{EQ-statement}]
Let $m\in \Nbb$ and $P\se \binom{[m]}{2}$ and $s\geq m^3$. We are required to show that there is an $([m],s)$-housing sequence of length at most $27$ and size at most $400\,m^4$ whose allocation set is $EQ(P)$. 

We make the assumption that $m$ is even and we omit the odd case since it is an easy adaptation of the even one. Our strategy will be to make $m^2$ copies for each of the $m$ variables (so $m^3$ in total), one for each possible pair $(i,j)\in[m]\times[m]$. Then for each $(i,j)\in P$, we will use one of the copies of the set of $m$ variables to perform the equalisation between variable $i$ and variable $j$. This way we will avoid crossing pairs.

Consider the following three sets of pairs in $\binom{[m^3]}{2}$ (see \autoref{fig:arching} for a concrete example):\vspace{-0.25cm}
\begin{align*}
    &R_1:=\;\{\{\,(2i-2)m+j,\,(2i-2)m+(2m-j+1)\,\}\;:\; i\in\{1,\ldots,m^2/2\}\;\text{and}\;j\in[m]\}\\
    &R_2:=\;\{\{\,(2i-1)m+j,\,(2i-1)m+(2m-j+1)\,\}\;:\; i\in\{1,\ldots,(m^2-2)/2\}\;\text{and}\; j\in[m]\}
\end{align*}

\vspace{-1cm}
\begin{align*}
    &R_3:=R_3'\cup R_3''\quad\text{where}\\
    &R_3':=\{\{\,xm^2+ym+x,\,xm^2+ym+y\,\}\,:\, \{x+1,y+1\}\in P,\;mx+y\;\text{even}\}\quad\text{and}\\
    &R_3'':=\{\{\,xm^2+ym+(m+1-x),\,xm^2+ym+(m+1-y)\,\}\,:\, \{x+1,y+1\}\in P,\;mx+y\;\text{odd}\}.\vspace{-0.4cm}
\end{align*}
Observe that $R_1$, $R_2$ and $R_3$ are each sets of disjoint pairs so that no two pairs in one of the sets cross. By \autoref{EQ-statement-2}, there are corresponding $([m^3],s)$-housing sequences $\Rcal_1$, $\Rcal_2$ and $\Rcal_3$ each of length at most $9$ and size at most $100 m^4$ with allocation sets $EQ(R_1),EQ(R_2)$ and $EQ(R_3)$, respectively. By \autoref{joint-allocation}, we get that the concatenation $\Rcal_1\Rcal_2\Rcal_3$ is of length $25$ and has allocation set $EQ(R_1)\cap EQ(R_2)\cap EQ(R_3)$.
We obtain our required housing sequence by appending a copy of $V(s,[m])$ to both the start and the end of $\Rcal_1\Rcal_2\Rcal_3$, thereby restricting our attention to the first $m$ of $m^3$ possible positions for houses. This $([m],s)$-housing sequence has allocation set $EQ(P)$ and length $27$, as required.
\end{pf}

\begin{pf}[Proof of \autoref{NEQ-statement}]
Similar as \autoref{EQ-statement-2} we prove the same statement with \lq $EQ(P)$\rq\ replaced by \lq $NEQ(P)$\rq. Using this statement in place of \lq \autoref{EQ-statement-2}\rq, we argue as in the proof of \autoref{EQ-statement} to deduce \autoref{NEQ-statement}.
\end{pf}

\begin{lem}\label{OR-statement-2}
For every $m\in 3\Nbb$, $s\geq m$ and $K\se[m]^3$ of triples of the form $(3f-2,3f-1,3f)$ for some $f\in[m/3]\}$, there is an $([m],s)$-housing sequence of length at most $15$ and size at most $300m^2$ with allocation set~$OR(K)$.
\end{lem}
\begin{pf}
    analogous to the proof of \autoref{EQ-statement-2} but using \autoref{OR-statement-special} instead of \autoref{EQ-statement-special}.
\end{pf}

\begin{pf}[Proof of \autoref{OR-statement}]
Let $m\in \Nbb$, $s\geq 27m^9$ and $K\se [m]^3$. We are required to show that there is an $([m],s)$-housing sequence of length at most $43$ and size at most $2001s^{2}$ whose allocation set is $OR(K)$.

Let $m'=m+3m^3$. For each $(i,j,k)\in[m]^3$, let $f((i,j,k))=m^2(i-1)+m(j-1)+(k-1)+m+1$. Consider the following set of (unordered) pairs from $\binom{[m+3m^3]}{2}$:\vspace{-0.2cm}
\[
R:=\{(i,3f),(j,3f+1),(k,3f+2):(i,j,k)\in[m]^3,f=f((i,j,k))\}\vspace{-0.2cm}
\]
That is, we pair each member of each possible triple $(i,j,k)\in[m]^3$ with the members of a unique consecutive triple $(3f,3f+1,3f+2)$ in the range $(m+1,\ldots,3m^3+m)$. By \autoref{EQ-statement}, there exists am $([m+3m^3],s)$-housing sequence $\Rcal$ whose length is at most $27$ and whose size is at most $4000 m^{12}$ and whose allocation set is $EQ(R)$. We think of $\Rcal$ as reserving for each $(i,j,k)\in[m]^3$ a consecutive triple somewhere in the range $(m+1,\ldots,3m^3+m)$.

Let $K'$ denote the set of triples $(3f,3f+1,3f+2)$ where $f=f(i,j,k)$ for each $(i,j,k)\in K$. By \autoref{OR-statement-2}, there exists a $([m+3m^3],s)$-housing sequence $\Kcal$ of length at most $15$ and size at most $1000m^{18}$ whose allocation set is $OR(K')$.

Then the concatenation $\Rcal\Kcal$ is a $([m+3m^3],s)$-housing sequence of length at most $41$ and size at most $2000m^{18}$ whose allocation set is $EQ(R)\cap OR(K')$. Our required housing sequence is obtained from $\Rcal\Kcal$ by appending the village $V([m],s)$ at either end. This housing sequence has index set $[m]$ and allocation set $OR(K)$.
\end{pf}

\section{Concluding remarks}
\label{sec:conc}

The main contribution of this paper is to establish NP-hardness for weak temporal sequences of 2-connected graphs. In the broader context of temporal graph minors, this result complements our earlier work~\cite{temporalsequences2025}. There, we introduced \emph{simultaneous embeddability} as a temporal analogue of planarity, but only after making certain restrictions on the admissible minor operations. The present work shows that those restrictions were in fact essential: if one relaxes them to allow weak temporal sequences, the corresponding embedding problem immediately becomes computationally intractable.

Thus the two canonical choices of graph minors suggested by the Graph Minor Structure Theorem lead to sharply different outcomes in the temporal setting. The \lq strict\rq\ variant yields a tractable theory~\cite{temporalsequences2025}, while the \lq weak\rq\ variant is already NP-hard (as proved here). Together with the extension to indefinite temporal sequences given in the appendix, this underlines that the notion of simultaneous embeddability used in~\cite{temporalsequences2025} is essentially the only one among these natural candidates that admits efficient algorithms.
For an outline of follow-up problems of this research we refer the reader to~\cite{temporalsequences2025}.

\newpage


\section{Appendix: Indefinite temporal sequences}
\label{sec:indef}
Let $\Gcal=(G_1,\ldots,G_n)$ be a (strict) temporal sequence. Then the \emph{indefinite temporal sequence} corresponding to $\Gcal$ is the sequence of graphs $(G_1,\ldots,G_n)$ without the data specifying which types of minor operations (deletion or contraction) are used to obtain the graphs from each other. An \emph{indefinite simultaneous embedding} of an indefinite temporal sequence $(G_1,\ldots,G_n)$ is a sequence of respective embeddings $(\iota_1,\ldots,\iota_n)$ so that whenever $G_i$ is a minor of $G_{i\pm1}$, there is a choice of $C,D\subseteq E(G_{i\pm1})$ so that $G_i=G_{i\pm1}/C\sm D$ and the embedding $\iota_i$ is the embedding induced by $\iota_{i\pm1}$ under these operations.
In this section, we outline a way to adapt the proof of \autoref{thm:main} to obtain the following.
\begin{prop}\label{thm:main_2}
    The problem of deciding whether indefinite temporal sequences of $2$-connected graphs admit weak simultaneous embeddings is NP-hard.
\end{prop}
For the proof outline, we will have to design appropriate analogues to the `Equaliser Gadget' (\autoref{fig:equaliser}), `Negator Gadget' (\autoref{fig:negator}) and `Or Gadget' (\autoref{fig:fror}). Since the first two such gadgets are strict temporal sequences, where the only minor operation is deletion, we can just reuse them. It is an elementary exercise to check that they behave analogously in the indefinite simultaneous embedding setting. In fact, aside from the `Or Gadget' and the machinery in \autoref{sec:overarching}, we can use identical argumentation to in \autoref{sec:part1}, \autoref{sec:part2} and \autoref{subsec:two}. An appropriate analogue for \autoref{sec:overarching} can also be proved; here we require that the pairs of sequences that we unify do not share edges where it is ambiguous as to what minor operations apply to them (instead of assuming their sets of weak edges sets are disjoint).
Hence, it remains to prove the following lemma which gives an analogue for the `Or Gadget'. A \emph{semi housing sequence} is obtained from a housing sequence by removing the final graph. The \emph{allocation set} of a semi housing sequence is defined analogously to the allocation set for a housing sequence.
\begin{lem}\label{lem:new_or}
    For every $m\in 3\Nbb$ and $i\in [m/3]$, there is an indefinite semi $(\{3i-2,3i-1,3i\},m)$-housing sequence of length $9$ and size at most $300\,m$ with allocation set $OR(\{3i-2,3i-1,3i\})$ and whose set of arches is $\{(3i-2,3i-1),(3i-1,3i)\}$.
\end{lem}
\begin{pf}
We will prove the lemma only in the case $(3i-2,3i-1,3i)=(1,2,3)$ since the other cases are analogous. The sequence $\Gcal$ is defined in \autoref{fig:amulet_1}.

\begin{figure}
    \begin{minipage}{0.45\textwidth}
        \large\scalebox{0.9}{\input{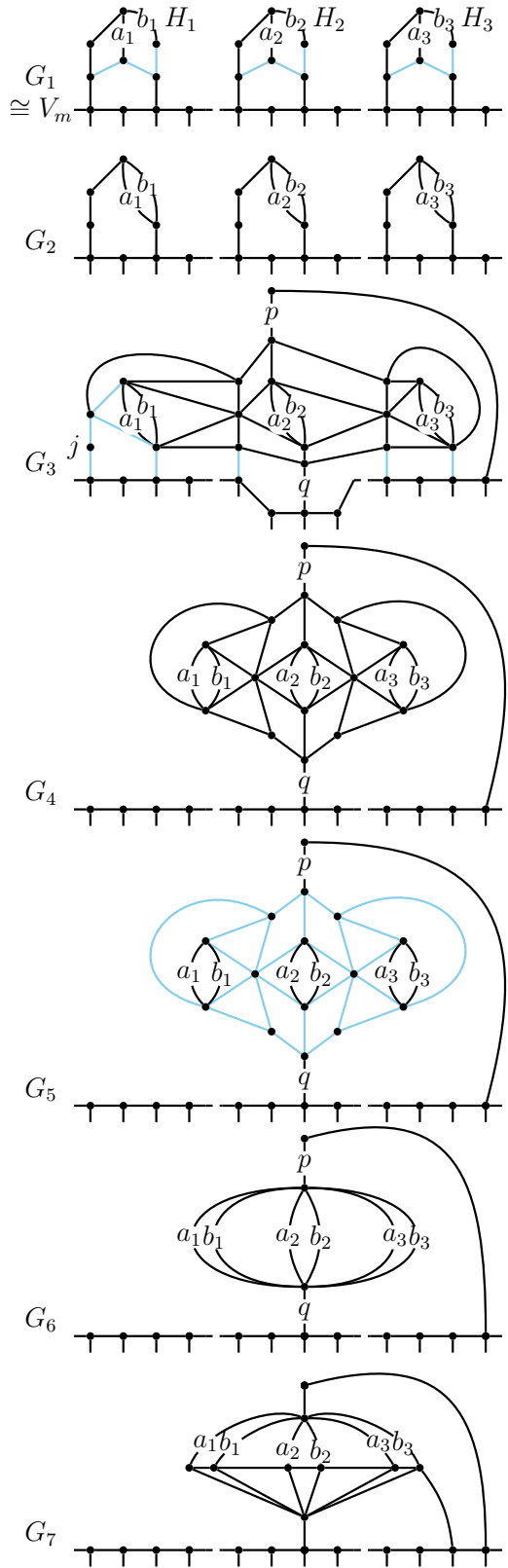}}
    \end{minipage}
    \hfill%
    \begin{minipage}{0.53\textwidth}
        \caption{The ``Indefinite Or Gadget'': an indefinite temporal sequence $\Gcal=(G_1\minorm\ldots\minorp G_7)$, analogous in some ways to a housing sequence. Illustrated are the houses $H_1$, $H_2$ and $H_3$ in the graph $G_1$, which is the village of equator $m$ and three houses. In $G_1$, the edges $a_1,b_1,a_2,b_2,a_3,b_3$ have been marked so that it is easier to keep track of the operations in this sequence. The graph $G_2$ is obtained from $G_1$ by minoring the \bcol\ edges; one may observe that there is a unique way to make the deletion/contraction choice that achieves this. The graph $G_3$ is obtained from $G_2$ by adding some vertices, adding some edges and uncontracting the edge~$q$. The graph $G_4$ is obtained from $G_3$ by minoring the \bcol\ edges; observe the deletion/contraction choice which achieves this is unique up to the operations on the edges incident with the vertex $j$. The graph $G_5$ is identical to $G_4$. The graph $G_6$ is obtained from the graph $G_5$ by minoring the \bcol\ edges; in this case, there are multiple different ways to choose the edges to delete/contract. The graph $G_7$ is obtained from $G_6$ by un-minoring some edges; one may observe that there is a unique way to make the deletion/contraction choice here as well. 
        Below is depicted an enlarged and recoloured version of the subgraph $A$ of $G_5$ ($=G_4$) induced on the \bcol\ edges (in the smaller version on the left) and the edges labeled $a_1,b_1,a_2,b_2,a_3,b_3,p,q$, and an enlarged version of the subgraph $B$ of $G_6$ induced on the edges labeled $a_1,b_1,a_2,b_2,a_3,b_3,p,q$.}

        \vspace{1em}\hspace{-2em}
        {\large\scalebox{0.9}{\begin{tikzpicture}
	\begin{pgfonlayer}{nodelayer}
		\node [style=black dot] (1118) at (-9.25, -1) {};
		\node [style=black dot] (1119) at (-9.25, 1) {};
		\node [style=black dot] (1120) at (-7.5, 0) {};
		\node [style=black dot] (1121) at (-5.75, 1) {};
		\node [style=black dot] (1122) at (-5.75, -1) {};
		\node [style=black dot] (1123) at (-4, 0) {};
		\node [style=black dot] (1124) at (-2.25, 1) {};
		\node [style=black dot] (1125) at (-2.25, -1) {};
		\node [style=black dot] (1126) at (-7.25, 2.25) {};
		\node [style=black dot] (1127) at (-7.25, -2.25) {};
		\node [style=black dot] (1128) at (-4.25, 2.25) {};
		\node [style=black dot] (1129) at (-4.25, -2.25) {};
		\node [style=black dot] (1130) at (-5.75, 3.5) {};
		\node [style=black dot] (1131) at (-5.75, -3.5) {};
		\node [style=small back] (1132) at (-1.75, 0) {};
		\node [style=small back] (1133) at (-2.75, 0) {};
		\node [style=small back] (1134) at (-5.25, 0) {};
		\node [style=small back] (1135) at (-6.25, 0) {};
		\node [style=small back] (1136) at (-8.75, 0) {};
		\node [style=small back] (1137) at (-9.75, 0) {};
		\node [style=none] (1138) at (-9.75, 0) {$a_1$};
		\node [style=none] (1139) at (-8.75, 0) {$b_1$};
		\node [style=none] (1140) at (-6.25, 0) {$a_2$};
		\node [style=none] (1141) at (-5.25, 0) {$b_2$};
		\node [style=none] (1142) at (-2.75, 0) {$a_3$};
		\node [style=none] (1143) at (-1.75, 0) {$b_3$};
		\node [style=none] (1144) at (-5.25, 3.75) {$x$};
		\node [style=none] (1145) at (-5.25, -3.75) {$y$};
		\node [style=none] (1146) at (-9.5, 1.5) {$z_1$};
		\node [style=none] (1147) at (-6.25, 1.25) {$z_2$};
		\node [style=none] (1148) at (-2, 1.5) {$z_3$};
		\node [style=none] (1149) at (-9.5, -1.5) {$w_1$};
		\node [style=none] (1150) at (-6.25, -1.25) {$w_2$};
		\node [style=none] (1151) at (-2, -1.5) {$w_3$};
		\node [style=none] (1152) at (-7.75, -0.5) {$u$};
		\node [style=none] (1153) at (-3.75, -0.5) {$v$};
		\node [style=none] (1154) at (-1.75, -3) {$A$};
		\node [style=black dot] (1155) at (-5.75, 5) {};
		\node [style=black dot] (1156) at (-5.75, -5) {};
		\node [style=small back] (1157) at (-5.75, 4.25) {};
		\node [style=none] (1158) at (-5.75, 4.25) {$p$};
		\node [style=small back] (1159) at (-5.75, -4.25) {};
		\node [style=none] (1160) at (-5.75, -4.25) {$q$};
		\node [style=black dot] (1161) at (4.25, 2) {};
		\node [style=black dot] (1162) at (4.25, -2) {};
		\node [style=black dot] (1163) at (4.25, 3.5) {};
		\node [style=black dot] (1164) at (4.25, -3.5) {};
		\node [style=small back] (1165) at (7.25, 0) {};
		\node [style=small back] (1166) at (6.25, 0) {};
		\node [style=none] (1167) at (6.25, 0) {$a_3$};
		\node [style=none] (1168) at (7.25, 0) {$b_3$};
		\node [style=small back] (1169) at (4.75, 0) {};
		\node [style=small back] (1170) at (3.75, 0) {};
		\node [style=none] (1171) at (3.75, 0) {$a_2$};
		\node [style=none] (1172) at (4.75, 0) {$b_2$};
		\node [style=small back] (1173) at (2.25, 0) {};
		\node [style=small back] (1174) at (1.25, 0) {};
		\node [style=none] (1175) at (1.25, 0) {$a_1$};
		\node [style=none] (1176) at (2.25, 0) {$b_1$};
		\node [style=small back] (1177) at (4.25, 2.75) {};
		\node [style=none] (1178) at (4.25, 2.75) {$p$};
		\node [style=small back] (1179) at (4.25, -2.75) {};
		\node [style=none] (1180) at (4.25, -2.75) {$q$};
		\node [style=none] (1183) at (7.25, -3) {$B$};
		\node [style=none] (1184) at (-7.25, 2.75) {$r$};
		\node [style=none] (1185) at (-1, 0.75) {};
	\end{pgfonlayer}
	\begin{pgfonlayer}{edgelayer}
		\draw [style=c0] (1119) to (1120);
		\draw [style=c0] (1120) to (1122);
		\draw [style=c0] (1122) to (1123);
		\draw [style=c0] (1123) to (1124);
		\draw [style=c0] (1123) to (1125);
		\draw [style=c2] (1123) to (1121);
		\draw [style=c0] (1121) to (1120);
		\draw [style=c0] (1120) to (1118);
		\draw [style=c0] (1130) to (1121);
		\draw [style=c0] (1126) to (1130);
		\draw [style=c0] (1130) to (1128);
		\draw [style=c0] (1127) to (1131);
		\draw [style=c0] (1131) to (1122);
		\draw [style=c2] (1131) to (1129);
		\draw [style=c0, in=165, out=150, looseness=2.50] (1126) to (1118);
		\draw [style=c0, in=15, out=30, looseness=2.75] (1128) to (1125);
		\draw [style=c0] (1118) to (1127);
		\draw [style=c0] (1129) to (1125);
		\draw [style=c0] (1120) to (1127);
		\draw [style=c2] (1129) to (1123);
		\draw [style=c0] (1123) to (1128);
		\draw [style=c0] (1126) to (1120);
		\draw [style=c0, bend left=45] (1119) to (1118);
		\draw [style=c0, bend left=45] (1118) to (1119);
		\draw [style=c0, bend left=45] (1121) to (1122);
		\draw [style=c0, bend left=45] (1122) to (1121);
		\draw [style=c0, bend left=45] (1124) to (1125);
		\draw [style=c0, bend left=45] (1125) to (1124);
		\draw [style=c0] (1119) to (1126);
		\draw [style=c0] (1128) to (1124);
		\draw [style=c0] (1155) to (1130);
		\draw [style=c0] (1131) to (1156);
		\draw [style=c0] (1163) to (1161);
		\draw [style=c0, bend left] (1161) to (1162);
		\draw [style=c0] (1162) to (1164);
		\draw [style=c0, bend left] (1162) to (1161);
		\draw [style=c0, bend left=75, looseness=1.75] (1161) to (1162);
		\draw [style=c0, bend left=75, looseness=1.75] (1162) to (1161);
		\draw [style=c0, bend left=90, looseness=2.50] (1161) to (1162);
		\draw [style=c0, bend left=90, looseness=2.50] (1162) to (1161);
	\end{pgfonlayer}
\end{tikzpicture}}}
    \label{fig:amulet_1}
    \end{minipage}
\end{figure}

Let $\Hcal$ be a subsequence of $\Gcal$ that contains the graph $G_5$. We say that an indefinite simultaneous embedding of $\Hcal$ is a \emph{right simultaneous embedding} if the embedding induced on the planet is the one depicted in \autoref{fig:village}. That is, from every pair of simultaneous embedding and its reorientation, we pick one. It is enough just to work with the right simultaneous embeddings of $\Gcal$, since we only care about allocation sets.
We first observe that the subsequence $(G_1\minorm G_2\minorp G_3\minorm G_4\minorp G_5)$ has exactly eight possible right simultaneous embeddings. Each right simultaneous embedding $(\sigma_1,\sigma_2,\sigma_3,\sigma_4,\sigma_5)$ is characterised by the triple $(T_{\sigma_1}(H_1),T_{\sigma_1}(H_2),T_{\sigma_1}(H_3))$. One may also characterise the right simultaneous embeddings by whether $a_i$ comes before or after $b_i$ in the rotator at $z_i$ in $G_4=G_5$, for each $i\in\{1,2,3\}$, in the embedding $\sigma_4=\sigma_5$.

For an embedding $\alpha$ of $A$, we say that it is \emph{right} if it is induced by an embedding of $G_5$ that comes from a right simultaneous embedding of the subsequence $(G_1\minorm G_2\minorp G_3\minorm G_4\minorp G_5)$.
We hence identify the right embeddings of $A$ with with the triples $(Q_1,Q_2,Q_3)\in\{0,1\}^3$.

Denote by $E'$ the set of edges in $A$ but not in $B$. That is, the \bcol\ ones in $G_5$ in \autoref{fig:amulet_1}. Let $\beta$ denote the embedding of $B$ illustrated in \autoref{fig:amulet_1}. We say that an embedding $\sigma$ of $A$ is \emph{admissible} if there is a choice of disjoint sets $D$ and $C$ so that $A\setminus D\;/\;C=B$ and the embedding induced on $B$ by $\sigma$ under this choice of operations yields $\beta$.

\begin{claim}\label{lem:admissible}
    The embeddings of $A$ induced by simultaneous embeddings of the subsequence triple $(G_5\minorm G_6\minorp G_7)$ are exactly the admissible embeddings of $A$.
\end{claim}
\begin{cproof}
    Observe that since it is $3$-connected, the graph $G_7$ admits a unique embedding $\beta'$ up to reorientation. There is a unique set of minor operations from which one may obtain $G_6$ from $G_7$. Furthermore, the embedding induced on $B$ under these operations by $\beta'$ is $\beta$.
\end{cproof}

\begin{figure}
    \large\scalebox{0.75}{\input{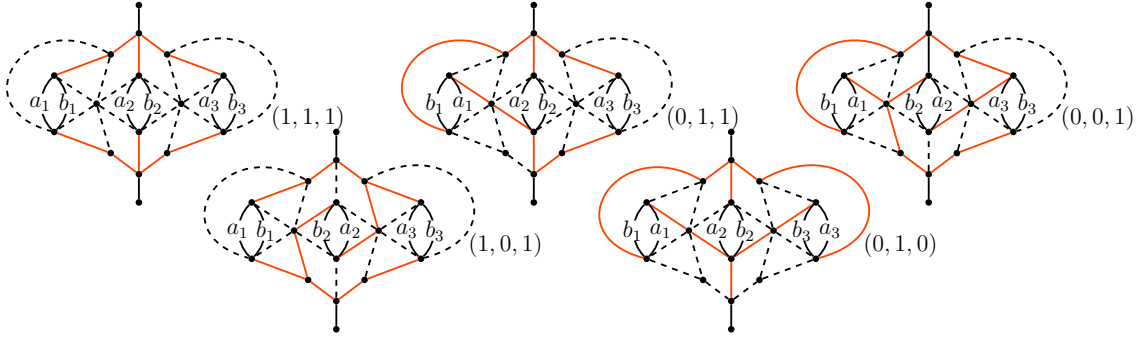}}
    \caption{A depiction of some admissible right embeddings of $A$. Depicted are the right embeddings associated with the triples $(1,1,1)$, $(0,1,1)$, $(0,0,1)$, $(1,0,1)$ and $(0,1,0)$. The precise choice of edges which must be deleted/contracted to produce the required embedding of $B$ are also marked: dashed edges should be deleted and edges coloured \rcol\ should be contracted. By symmetry, we also have that the embeddings associated with the triples $(1,1,0)$ and $(1,0,0)$ are admissible. This shows that all triples are admissible except possibly $(0,0,0)$.}
    \label{fig:amulet_2}
    \vspace{-0.3cm}
\end{figure}

In \autoref{fig:amulet_2} we prove that all right embeddings except for the one associated with $(0,0,0)$ are admissible. The following lemma completes the characterisation.

\begin{claim}\label{lem:not_admissible}
    The right embedding of $A$ associated with the triple $(0,0,0)$ is not admissible.
\end{claim}
\begin{cproof}
    We will make use of the labeling of the graph $A$ given in \autoref{fig:amulet_1}. Let $\sigma$ denote the embedding of $A$ corresponding to the triple $(0,0,0)$. Let $D$ and $C$ be disjoint sets so that $D\sqcup C=E'$. Further suppose that $A/C\setminus D=B$ (up to relabeling vertices). Let $\beta'$ denote the embedding induced by $\sigma$ under the given operations. We are required to show~$\beta'\neq\beta$. 

    Let $x'$ and $y'$ be the vertices in $A$ whose branch sets contain $x$ and $y$, respectively, and $T_x$ and $T_y$, respectively, be the disjoint subgraphs of $A$ induced on these branch sets. Now suppose for a contradiction, that for each $i\in\{1,2,3\}$, the rotator $\beta'(x')$ contains as a restriction the cyclic orientation $pa_ib_i$ (and not $pb_ia_i$). Then the subgraph $T_x$ contains $w_1,w_2$ and $w_3$ and the subgraph $T_y$ contains $z_1,z_2$ and $z_3$. Hence there is a path $Q$ in $T_y$ from $y$ to $z_2$ and it uses either $p$ or $q$, without loss of generality assume $q$. The path $Q$ is marked by colouring its edges \rcol\ in \autoref{fig:amulet_1}. Now observe that there must be a path $P$ in $T_x$ between $x$ and $w_2$. Observe that $P$ must use both $p$ and $r$. But then there is no path from $z_1$ to $y$ that avoids $P$, a contradiction.
\end{cproof}

    Together, \ref{lem:not_admissible} and \ref{lem:admissible} prove the allocation set for $\Gcal$ to be $OR(\{1,2,3\})$, as required.
\end{pf}

\bibliographystyle{plain}
\bibliography{wliterature}

\end{document}